\documentclass[preprint,12pt]{elsarticle}

\usepackage{amssymb}
\usepackage{amsmath}
\usepackage{amsthm}

\usepackage{lineno}

\usepackage{physics}
\usepackage{siunitx}
\usepackage{subcaption}
\usepackage[fleqn]{mathtools}
\usepackage{multirow}
\usepackage{booktabs}
\usepackage{url}
\usepackage{algorithm}
\usepackage{algpseudocode}

\usepackage{tikz}
\usepackage{eso-pic}
\usepackage{xcolor}

\algrenewcommand\algorithmicrequire{\textbf{Inputs:}}
\algrenewcommand\algorithmicensure{\textbf{Outputs:}}

\journal{Acta Astronautica}

\renewcommand\fbox{\fcolorbox{black}{gray!10}}

\newcommand\acceptedmanuscripttext{%
	\normalfont\scriptsize
	\textbf{Accepted Manuscript.} 
	This manuscript has been accepted for publication in
	\textit{Acta Astronautica}.\\
	It has not undergone final copy-editing, typesetting, or proofreading and may be outdated.\\[5pt]
	The final Version of Record is available at
	\url{https://doi.org/10.1016/j.actaastro.2026.07.020} and is distributed under the terms of a CC BY-NC-ND 4.0 license.\\[5pt]
	Please cite the published version as:\\
	\textbf{A. Nunes, S. Brás, P. Batista, J. Xavier, Designing trajectories in the Earth--Moon system: A Levenberg--Marquardt approach, Acta Astronautica 249 (2026) 337--357.}
}

\AddToShipoutPicture*{%
	\AtPageUpperLeft{%
		\begin{tikzpicture}[remember picture,overlay]
			\node[
			anchor=north,
			yshift=-10pt
			] at (current page.north)
			{%
				\fbox{%
					\parbox{\dimexpr\textwidth-2\fboxsep-2\fboxrule}{%
						\acceptedmanuscripttext
					}%
				}%
			};
		\end{tikzpicture}
	}%
}

\begin{document}

\begin{frontmatter}


\title{Designing trajectories in the Earth--Moon system:\\ a Levenberg--Marquardt approach}

\author[ISR]{António Nunes\corref{ANcor}} 
\ead{antonio.w.nunes@tecnico.ulisboa.pt}
\cortext[ANcor]{Corresponding author at: Institute for Systems and Robotics, Instituto Superior Técnico, Universidade de Lisboa, Portugal}

\author[ESA]{Sérgio Brás} 
\author[ISR]{Pedro Batista} 
\author[ISR]{João Xavier} 

\affiliation[ISR]{organization={Institute for Systems and Robotics, Instituto Superior Técnico, Universidade de Lisboa},
	addressline={Av. Rovisco Pais 1}, 
	postcode={1049-001 },
	city={Lisbon},
	country={Portugal}}

\affiliation[ESA]{organization={European Space Agency},
	addressline={Keplerlaan 1, PO Box 299}, 
	postcode={2200AG},
	city={Noordwijk},
	country={The Netherlands}}

\begin{abstract}
Trajectory design in cislunar space under a high-fidelity ephemeris model (HFEM) is pursued through a nonlinear optimization perspective anchored on the transition of solutions from lower fidelity models, namely the Circular Restricted Three-Body Problem (CR3BP). The optimization problem is posed in the likeness of a multiple-shooting approach, aiming for segment-to-segment continuity while tracking proximity to the original CR3BP structures. The analysis of various formulations leads to the selection of an unconstrained least-squares problem for further investigation. The nonlinear optimization problem is convexified and the use of the Levenberg--Marquardt algorithm, as an alternative to the minimum-norm update equation found in most literature, is investigated for its control over the update step and inherent robustness. Additional techniques, such as adaptive weighting, are employed to further consolidate the behavior of the proposed algorithm in challenging scenarios. Numerical trials evaluate the adequacy of the methodology presented and compare it to the minimum-norm baseline over various application cases, including the generation of quasi-periodic trajectories and orbital transfers between them. The proposed technique is found to be a suitable alternative to the minimum-norm scheme, generally retaining better proximity to the original CR3BP trajectories and providing benefits in numerical robustness and stability. Moreover, the ease of including proximity objectives in a relaxed manner is shown to facilitate control over the shape of the final converged solution.
\end{abstract}


\end{frontmatter}


\section{Introduction}

Consolidating a human presence in cislunar space is a topic that has returned to the main stage of scientific research on the sustainable exploration of space and its assets. Motivated by the recent efforts of international space agencies, such as NASA and ESA, returning mankind to the lunar surface has been given upmost priority with initiatives like the Artemis Program currently underway and proposed plans for Gateway, a lunar orbiting station \cite{creech2022ArtemisOverview,lehnhardt2024Gateway}. Evidently, the design of optimal trajectories to be followed by spacecraft within the Earth--Moon system is fundamental for ensuring that operational objectives may be met while respecting physical constraints and minimizing control efforts. In turn, this boosts mission longevity and feasibility, given how every extra kilogram of fuel translates into reduced payload capacity or a significant increase in cost.

Given the inherently chaotic nature of the dynamics at play, any modern trajectory design algorithm ought to consider high-fidelity ephemeris models (HFEMs) to accurately capture the complex relative motion of the Earth and Moon and prominent perturbing effects. Nonetheless, simpler approximate representations of the dynamics, such as the Circular Restricted Three-Body Problem (CR3BP), often provide relevant insights on the dynamical structures at play \cite{gomez2002SSFrequencies}, enabling a straightforward generation of periodic and quasi-periodic orbits that may be considered in the preliminary stages of mission design. Naturally, when more realistic dynamical representations are considered, the symmetries that make the CR3BP particularly tractable are destroyed, rendering the process of determining trajectories that meet the operational goals of a particular mission in a HFEM non-trivial. Despite this, it is oftentimes the case that high-fidelity dynamics can be interpreted as perturbations to the CR3BP, such that the more intricate solutions of HFEMs exist in the vicinity of CR3BP counterparts \cite{park2025CharacterizationL2Analogs, Soto2025Persistence}. This ultimately means that solutions from the circular problem typically provide convenient initial guesses for trajectory design algorithms, which can hence be interpreted as \textit{transitioning} these low-fidelity solutions into high-fidelity dynamical models.

The use of so-called \textit{shooting} methods dominates recent research on the generation of orbits in the CR3BP and their transition to other dynamical models. The simplest implementations correspond to \textit{single-shooting} methods, which consider a single point in phase space (i.e., in terms of spacecraft position and velocity) from which the dynamics are propagated. This point is successively updated until a desired trajectory is achieved over a given window of time, which may itself also be a variable of the design problem. The simplicity of single-shooting methods makes them especially suitable for the purposes of finding periodic orbits under the CR3BP, namely by exploiting dynamical symmetries in space according to the Mirror Theorem \cite{RoyMirror1955Theorem}, as first approached in \cite{howell1984HaloOrbits}. Equivalently, single-shooting methods operate well at transitioning CR3BP orbits to models of intermediate complexity, such as the Elliptic Restricted Three-Body Problem (ER3BP), where some of the dynamical symmetries are maintained and the solutions from the circular problem offer very adequate initial guesses \cite{ferrari2018OrbitsER3BP}.

However, while the ER3BP represents a step-up in realism in comparison to the CR3BP, it does not provide the high degree of accuracy that can only be achieved through a HFEM representation. In particular, high-fidelity models are characterized by non-autonomous dynamics with no lasting symmetries, which means they do not allow for the existence of purely periodic solutions. Hence, trajectories under these dynamics need to be fully characterized over the entire time span of a mission, which may amount to several months or even years \cite{williams2017NRHO,Lee2019NRHO,Spreen2023BaselineNRHO}. This concern is worsened by the fact that the already chaotic nature of the CR3BP is exacerbated in the time-varying dynamics of a HFEM, making it especially challenging to precisely integrate a trajectory from a single point over a large window of time, which renders single-shooting techniques impractical, from a numerical perspective, for most applications. With this in mind, trajectory design under HFEMs is most commonly pursued through \textit{multiple-shooting} methods that divide a trajectory into various segments. Since the initial CR3BP trajectory will not be feasible in the more realistic HFEM dynamics, separate integration of each segment will initially yield a discontinuous solution. Hence, the base objective of multiple-shooting algorithms is to guarantee continuity between segments, which is done by enacting slight changes to the position, velocity, and possibly the epoch of each \textit{patch point} --- typically the point at the start of each segment. For forward multiple-shooting methods, this means that the integration of each individual segment should be consistent with the succeeding segment  \cite{park2025CharacterizationL2Analogs,Spreen2023BaselineNRHO}. Equivalently, other variations such as forward-backward multiple-shooting methods propagate the dynamics from a central point along each segment and the continuity constraint is imposed at both ends instead, yielding benefits in convergence at the cost of a higher computational burden \cite{williams2017NRHO, Diane2017stationkeeping}. While most of the literature deals with the problem of transitioning CR3BP periodic orbits into quasi-periodic HFEM counterparts \cite{park2025CharacterizationL2Analogs,williams2017NRHO,Lee2019NRHO,Spreen2023BaselineNRHO,Spreen2020NRHO}, work on transfer trajectories to, from, and between orbits in HFEMs is also available, evidencing the versatility of multiple-shooting algorithms \cite{Diane2017stationkeeping,lantoine2011transfers, pavlak2010thesis, Scheuerle2025LunarTransfers}. In addition, multiple-shooting has also been considered as a means of interconnecting trajectory planning with spacecraft station-keeping, allowing for discontinuities in velocity at specific locations along an orbit that are interpreted as impulsive control maneuvers, which may or may not be optimized \cite{Diane2017stationkeeping,pavlak2012Strategy,Guzzetti2017stationkeeping}.

Since there are typically more variables than constraints in a multiple-shooting approach, one expects various (possibly even infinite) solutions to be available. Hence, if the problem is not further constrained, there are usually no guarantees to which specific solution an algorithm will converge, and it may be the case that the converged trajectory is significantly far away from the initial CR3BP guess, which would most likely have been carefully chosen during a preliminary study. To possibly mitigate this effect, most approaches in the literature employ a minimum-norm (MN) scheme that, at each iteration, minimizes the norm of the update step to be applied to the variables of the optimization problem, with the intent of remaining close to the initial guess \cite{park2025CharacterizationL2Analogs,Belinda2007ImprovedCorrections}. Practically speaking, it is observed that this algorithm is extremely well suited over a wide range of applications, being commonly chosen over other alternatives due to its convergence speed and computational efficiency \cite{park2025CharacterizationL2Analogs}. However, we show that the minimization of the update step does not guarantee that proximity to the CR3BP orbit is maintained and that it may be non-trivial to leverage additional objectives in a general and robust manner. Moreover, the step-minimization approach pursued in most MN methods may not be sufficient to ensure the algorithm does not diverge. Augmented Lagrangian formulations \cite{DeiTos2017Refinement} or the use of \textit{line search} and \textit{trust region} techniques \cite{mangialardo2018tesi} have been presented as means of addressing this concern, though these are known to significantly increase the computational burden of the application. Alternatively, some authors propose manually limiting the norm of the update step according to a maximum threshold \cite{Spreen2023BaselineNRHO}, whose best value might however not be evident to predict for each application. Furthermore, this procedure does not fully mitigate the possibility of numerical divergence that may affect the MN scheme.

Current research is also particularly interested in leveraging dynamical insights onto the design of trajectories within HFEMs, e.g. by exploiting manifold information from simpler models such as the CR3BP \cite{Spreen2020NRHO}. Alternatively, homotopy-based approaches leveraging dynamical continuation of incrementally realistic models, for example the ER3BP \cite{park2025CharacterizationL2Analogs, park2022IntermediateDynamicalModels,Gao2026Continuation} or different versions of the Restricted Four-Body Problem \cite{park2022IntermediateDynamicalModels,boudad2020BC4BP, park2025bridgingHFEM,Sanaga2025Leveraging}, have also become prevalent in the literature. These techniques essentially bridge low- and high-fidelity representations of the dynamics, easing the complexity of transitioning solutions from one to the other. In particular, recent studies such as the Unified Transition Scheme presented in \cite{sanaga2025UTS,sanaga2026UTSTori}, establish a coherent framework connecting models of varying fidelity, investigating different continuation parameters and the incorporation of underlying structures in the transition process, such as invariant tori \cite{sanaga2026UTSTori,Brown2025InvariantTori}. Frequency analysis as a way of interpreting HFEM dynamics as a perturbation to the CR3BP, to exploit conservation of major dynamical frequencies in the orbit transition process, is also a topic of recent research \cite{park2025CharacterizationL2Analogs,DeiTos2017Refinement,park2025Frequency}, namely based on previous work to identify major perturbing frequencies in the Solar System \cite{gomez2002SSFrequencies}. Assessing the fundamental structures behind each dynamical representation has also been of utmost importance for understanding the cases where transition is not possible, namely through the viewpoint of bifurcation theory \cite{park2025CharacterizationL2Analogs,Gao2026Continuation,spreen2021phd}.

Most of the research has, however, successively strayed away from an optimization problem viewpoint. In particular, the underlying scheme used to solve the nonlinear multiple-shooting problem behind a typical transition problem has not seen significant developments beyond the original implementations of the MN update step. Notable exceptions include the use of nonlinear optimization solvers, with some works repurposing algorithms developed for previous missions, such as Copernicus \cite{williams2017NRHO,williams2012copernicus}, and others considering more modern alternatives like SNOPT \cite{gill2002SNOPT}, a sparse solver that exploits the sparsity of the problem's Jacobian, as approached in \cite{oguri2020EQUULEUS}. Nevertheless, these solvers are typically rather complex to implement and computationally demanding to run. In turn, this motivates the development of alternative numerical schemes that retain the implementation simplicity and computational lightness of a typical MN update step, while striving for improved numerical robustness and flexibility at handling multiple objectives.

In light of the aforementioned considerations, we return to the design of the nonlinear optimization problem at hand, building upon the seminal work of MN schemes from literature to cover different formulations and introduce a generic proximity constraint to control the shape of the final solution. This study leads to a novel algorithm for HFEM transition that is lightweight, flexible in accommodating mission objectives, and robust to divergence, thus addressing the drawbacks identified in the existing literature. In particular, an update step damping term is introduced to guarantee a monotonous, non-increasing evolution of the continuity residual, pursuing a Levenberg--Marquardt (LM) line of reasoning \cite{Levenberg1944,Marquardt1963} --- a key algorithm in the field of nonlinear least-squares estimation, with applications to the assessment of near-Earth spacecraft orbits and attitude from available measurements \cite{warner2016OrbitDetermination,alarcon2005OrbitDetermination,scire2015OrbitDetermination,wu2025AttitudeEstimation}. The proposed algorithm maintains an analytical update equation, ensuring a general and straightforward application. Its adequacy over a wide range of test cases is carefully analyzed and compared with the MN alternative presented in \cite{Spreen2023BaselineNRHO}, which serves as baseline in terms of convergence speed and proximity to the initial guess. Additional techniques, such as adaptive weighting, are leveraged to improve convergence of the proposed methodology when multiple objectives are imposed.  In sum, the proposed LM algorithm is assessed as a potential alternative to the MN update step for its inherent robustness and versatility, which may be leveraged by practitioners when addressing challenging application cases.

The study carried out in this work is structured as follows. Firstly, Section~\ref{sec:dynamical_models} provides the theoretical background and analytical description of the dynamical models at play, i.e. the CR3BP and HFEM, highlighting key differences and possible challenges in the initialization of the trajectory transition process from CR3BP solutions. In Section~\ref{sec:trajectory_design}, the trajectory design problem is discussed and formulated via a multiple-shooting approach, drawing attention to relevant constraints to be imposed. Section~\ref{sec:op_probs} builds on these concepts, framing multiple-shooting as a nonlinear optimization problem under different formulations of varying constrained nature. These are then convexified according to the description provided in Section~\ref{sec:convex}, leading to the MN update equation used in most of the literature and an alternative approach based on the LM algorithm. Relevant extensions to the proposed implementation are also discussed. The complete algorithm and initialization procedure are detailed and schematically laid out in Section~\ref{sec:algorithm}. In Section~\ref{sec:results}, the proposed approach is evaluated under various application scenarios concerning the generation of quasi-periodic trajectories and orbital transfers between them. A comparison analysis is performed with respect to the baseline MN update equation. Finally, Section~\ref{sec:conclusions} provides a discussion on the results obtained and the most meaningful conclusions that may be drawn from this work.

\section{Dynamical models}
\label{sec:dynamical_models}
As previously introduced, the trajectory design process in cislunar space is typically based on the transition of solutions from low-fidelity dynamical models to more realistic, albeit more complex, representations of the dynamics. In fact, since HFEMs may oftentimes be interpreted as perturbations to simpler models, such as the CR3BP \cite{gomez2002SSFrequencies}, it is known that solutions in realistic dynamics \textit{tend} to exist in the vicinity of those from low-fidelity models, provided that relevant Hamiltonian structures are preserved \cite{Soto2025Persistence}. Perhaps the most prominent and well-studied example is that of quasi-periodic trajectories under HFEMs that describe multiple windings in close proximity to truly periodic counterparts from the CR3BP. In \cite{park2025CharacterizationL2Analogs,Sanaga2024ChallengingRegion}, notable exceptions are studied for the case of Halo orbits around the L2 equilibrium point of the Earth--Moon system, identifying transition-challenging regions of space and highlighting difficulties at transitioning orbits whose period is close to resonance with the Earth--Moon dynamics, where trajectories tend to unfold into multiple lobes. However, for most general purposes, periodic orbits from the CR3BP offer sufficiently adequate initial guesses for the purposes of determining quasi-periodic counterparts under HFEMs, with some techniques being available to circumvent difficulties in challenging application scenarios. Similarly, it is typically possible to transition a generic solution from the circular problem into a higher-fidelity model. In this sense, it is relevant to briefly explore the differential equations that dictate spacecraft motion within the CR3BP and HFEM, to highlight key differences in the dynamical representations and types of solutions that may be found.

\subsection{Circular restricted three-body problem (CR3BP)}
\label{sec:CR3BP}

A restricted three-body problem considers, in its most lenient definition, the motion of a small mass (e.g., a spacecraft) subject to the gravitational forces of two much more massive bodies in its proximity, named primaries. Given the extreme differences in mass between primaries and spacecraft, the former move independently to the latter and thus constitute an isolated two-body problem, which is known to have a limited set of solutions \cite{szebehely1967TheoryOfOrbits}. Of particular interest is the case where the two bodies describe Keplerian orbits around their shared barycenter, since this provides an adequate first approximation for the motion of many two-body configurations within the Solar System, most notably the Earth--Moon binary studied in this work. The CR3BP corresponds to the specific case where the orbits of the primaries are further approximated to being circular. While this is clearly a step backwards from the pursuit of realism as compared to the ER3BP, for example, the CR3BP is widely esteemed from an academic standpoint since most two-body systems of interest are characterized by low orbital eccentricity (e.g., ${e\approx0.0549}$ for the Earth--Moon case) and because this model provides tractable but relevant insights on the dynamical structures at play.

One of the particularities of the CR3BP is that its equations of motion (EoM) are fully autonomous under the appropriate reference frame, which is made to rotate at a constant angular velocity matching the orbits of the primaries, hence being referred to as the \textit{synodic} reference frame. Aligning both primaries along the $X$ direction, the angular velocity vector with $Z$, and completing the right-handed referential with the definition of $Y$ perpendicular to both, the synodic reference frame is scaled to bring the distance between primaries to unity. Additionally, a time scale is defined such that the orbital period of the two massive bodies is brought to $2\pi$, resulting in unit angular velocity. Denoting the most massive primary by subscript $1$ and the other by subscript $2$, without loss of generality, the masses of the primaries are scaled to $m_1=1-\mu$ and $m_2=\mu$, where $\mu$ is a problem-specific constant. This procedure fixes the position of the first and second primaries to ${\mathbf{r}_1=(-\mu,0,0)}$ and $\mathbf{r}_2=(1-\mu,0,0)$, respectively. In this work, the Earth--Moon system is considered, for which $\mu\approx0.01215$, and hence subscript $1$ denotes the Earth and subscript $2$ the Moon.

Under the proposed scaling and reference frame, it is possible to establish an effective potential of the CR3BP as a function of the spacecraft's position, $\mathbf{r}^T = \begin{bmatrix} x & y & z \end{bmatrix}$, given by
\begin{equation}
	\label{eq:CR3BP_eff_potential}
	U = \frac{x^2 + y^2}{2} + \frac{1-\mu}{\norm{\mathbf{r}-\mathbf{r}_1}} + \frac{\mu}{\norm{\mathbf{r}-\mathbf{r}_2}}.
\end{equation}
This potential is characterized by a Coriolis acceleration term, due to the rotation of the synodic frame, and two gravitational contributions, from the primaries considered. Adopting the dot notation to represent time derivatives, and assuming no sensor noise or external perturbations, the resulting EoM for the spacecraft are written as
\begin{equation}
	\label{eq:CR3BP_eom}
	\begin{aligned}
		\dot{\mathbf{r}} &= \mathbf{I_3} \mathbf{v},\\
		\dot{\mathbf{v}} &= \Omega \mathbf{v} + \nabla U(\mathbf{r}),
	\end{aligned}
	\qquad \text{with} \qquad \Omega=\begin{bmatrix}
		0 & 2 & 0\\
		-2 & 0 & 0\\
		0 & 0 & 0
	\end{bmatrix},
\end{equation}
where $\mathbf{v}^T=\begin{bmatrix} \dot{x} & \dot{y} & \dot{z}\end{bmatrix}$ is the spacecraft velocity, $\mathbf{I_3}$ is the $3\times3$ identity matrix, and $\nabla U(\mathbf{r})$ is the gradient of the effective potential in Eq.~\eqref{eq:CR3BP_eff_potential}, evaluated at $\mathbf{r}$. 

Note how Eq.~\eqref{eq:CR3BP_eom} is independent of time, in accordance with the previous remarks on the autonomous nature of the CR3BP EoM. Moreover, since the EoM are also Hamiltonian, they obey the conservation of a (single) integral of motion,
\begin{equation*}
	\mathcal{J} = (x^2+y^2)+2\left(\frac{1-\mu}{\norm{\mathbf{r}_1}} + \frac{\mu}{\norm{\mathbf{r}_2}}\right)- \norm{\mathbf{v}}^2,
\end{equation*}
termed the Jacobi integral \cite{ross3BPbook}, which represents a measure of the spacecraft's total energy per unit mass in the synodic reference frame. These characteristics are responsible for the existence of infinite periodic solutions in phase space, as first pointed out by Poincaré \cite{poincare1892MethodesNouvelles}, which evidently makes the CR3BP a great tool for the early steps of a trajectory design process. 

\subsection{High-fidelity ephemeris model (HFEM)}
In comparison, a HFEM considers the gravitational influence of up to $P$ celestial bodies, whose true motion is accurately depicted based on past and predicted ephemeris data. Naturally, as the number of celestial bodies is increased, the dynamics get increasingly more complex and realistic. However, for the case of the Earth--Moon system, it is typically the case that only the gravitational influence of the Sun is of particular relevance, with the remaining celestial bodies exerting nearly negligible perturbations in most applications \cite{Park2024AssessmentDynModels}. As such, the HFEM employed in this work considers only the gravitational effect of the Earth, Moon, and Sun, whose motion is depicted via planetary ephemeris data.

Since the use of planetary ephemeris introduces time-dependent effects that cannot be counteracted through the definition of an appropriate reference frame, it is common to establish the EoM of the spacecraft in an inertial frame relative to a so-called central body, denoted by subscript $q$. In this work, the Earth-centered J2000 frame, a widely used standard, was selected. In this sense, neglecting external perturbations and considering the gravitational contributions of each celestial body under analysis as point-masses, the spacecraft's motion under HFEM dynamics is dictated by
\begin{equation}
	\label{eq:HFEM_eom}
	\begin{aligned}
		\dot{\mathbf{r}} &= \mathbf{I_3} \mathbf{v},\\
		\dot{\mathbf{v}} &= -\mu_q\frac{\mathbf{r}}{\norm{\mathbf{r}}^3} + \sum_{\substack{j=1\\ j\neq q}}^{P} \mu_j \left(\frac{\mathbf{r}_j - \mathbf{r}}{\norm{\mathbf{r}_j - \mathbf{r}}^3} - \frac{\mathbf{r}_{j}}{\norm{\mathbf{r}_{j}}^3}\right)
	\end{aligned},
\end{equation}
where $\mathbf{r}_j$ is the position of body $j$ in the J2000 frame, at each instant, and $\mu_j$ its gravitational parameter --- not to be confused with the problem-specific constant of the CR3BP.
Note that Eq.~\eqref{eq:HFEM_eom} may also be scaled according to relevant quantities specific to the set of celestial bodies being considered, which may be time-varying or fixed, but it is not uncommon for off-the-shelf HFEM propagators to work directly with dimensional units. In comparison with the CR3BP EoM, Eq.~\eqref{eq:HFEM_eom} depends directly on time through $\mathbf{r}_j$ without the introduction of relevant symmetries. This has the direct consequence of inhibiting the existence of truly periodic orbits.

\subsection{Coordinate transformations}
\label{subsec:coord_transf}
Given that the two dynamical models approached in this work differ so significantly at the level of the reference frame under which they are established, it makes sense to briefly go over how one may transform relevant quantities between the barycentric synodic reference frame of the CR3BP and the J2000 inertial reference frame considered in the HFEM. This is crucial for the purposes of constructing the initial guess for the trajectory design algorithms. For the sake of consistency, assume that all quantities presented in this section have been scaled by the CR3BP characteristic length, time, and linear velocity factors --- the latter of which may be directly inferred from the former two.

To discern between quantities expressed in the J2000 inertial reference frame, centered at the Earth, and the synodic reference frame of the CR3BP, centered at the Earth--Moon barycenter, we may use the prescript notations $J$ and $C$, respectively. Then, the coordinate transformation from the barycentric synodic frame to the J2000 inertial frame may be described, at each instant, by
\begin{equation}
	\label{eq:CR3BP_J2000_transf}
	\prescript{J}{}{\mathbf{r}} = \prescript{J}{}{\mathbf{r}}_{b} + \frac{d_\text{HFEM}}{d_\text{CR3BP}} \mathbf{R} \prescript{C}{}{\mathbf{r}},
\end{equation}
where $d_\text{HFEM}$ and $d_\text{CR3BP}$ denote the distance between the Earth and Moon in the HFEM and CR3BP dynamics, respectively, $\mathbf{R}$ is the rotation matrix between the two reference frames, and $\prescript{J}{}{\mathbf{r}}_{b}$ is the position of the Earth--Moon barycenter relative to Earth. The multiplicative factor $d_\text{HFEM}/d_\text{CR3BP}$ accounts for the varying distance between the Earth and Moon, which is not reflected in the simpler CR3BP dynamics, since $d_\text{CR3BP}$ is a fixed value (equal to one, using the scaling previously mentioned). The position of the barycenter in the J2000 frame, $\prescript{J}{}{\mathbf{r}}_{b}$, is given by
\begin{equation*}
	\prescript{J}{}{\mathbf{r}}_{b} = \frac{\mu_M}{\mu_E + \mu_M}\prescript{J}{}{\mathbf{r}}_{M},
\end{equation*}
where $\mu_E$ and $\mu_M$ are the gravitational parameters of the Earth and Moon, respectively, and $\prescript{J}{}{\mathbf{r}}_{M}$ is the position of the Moon relative to Earth. Finally, the rotation matrix $\mathbf{R}$ can be promptly determined, at each instant, through a series of simple computations involving the normalized base directions of the CR3BP synodic reference frame, as measured from the J2000 inertial frame. Following the exposition detailed in Section~\ref{sec:CR3BP}, these may be established as
\begin{equation*}
	\mathbf{x}_0 = \frac{\prescript{J}{}{\mathbf{r}}_{M}}{\norm{\prescript{J}{}{\mathbf{r}}_{M}}}, \qquad
	\mathbf{z}_0 = \frac{\prescript{J}{}{\mathbf{r}}_{M}\times\prescript{J}{}{\mathbf{v}}_{M}}{\norm{\prescript{J}{}{\mathbf{r}}_{M}\times\prescript{J}{}{\mathbf{v}}_{M}}}, \qquad \text{and} \qquad
	\mathbf{y}_0 = -\mathbf{x}_0\times\mathbf{z}_0,
\end{equation*}
where $\prescript{J}{}{\mathbf{v}}_{M}$ denotes velocity of the Moon relative to Earth. With this information, the rotation matrix is given by
\begin{equation*}
	\mathbf{R} = \begin{bmatrix}
		\mathbf{x}_0 & \mathbf{y}_0 & \mathbf{z}_0
	\end{bmatrix}.
\end{equation*}

Differentiating Eq.~\eqref{eq:CR3BP_J2000_transf} leads to the transformation of spacecraft velocity components, i.e.
\begin{equation}
	\label{eq:CR3BP_J2000_transf_v}
	\prescript{J}{}{\mathbf{v}} = \prescript{J}{}{\mathbf{v}}_b + \left(\frac{\dot{d}_\text{HFEM}}{d_\text{CR3BP}}\mathbf{R} + \frac{d_\text{HFEM}}{d_\text{CR3BP}}\dot{\mathbf{R}}\right)\prescript{C}{}{\mathbf{r}} + \frac{d_\text{HFEM}}{d_\text{CR3BP}}\mathbf{R}\prescript{C}{}{\mathbf{v}},
\end{equation}
where $\prescript{J}{}{\mathbf{v}}_b$ is the velocity of the Earth--Moon barycenter relative to Earth. The instantaneous variations $\dot{d}_\text{HFEM}$ and $\dot{\mathbf{R}}$ may be computed numerically via finite differences, for example. Note that Eq.~\eqref{eq:CR3BP_J2000_transf_v} assumes a uniform correspondence of the passing of time between the CR3BP and HFEM. This is a topic deserving of more attention, and its discussion will follow shortly.

The inverse transformations, i.e. those from the CR3BP synodic reference frame to the J2000 inertial frame, may be established directly by reorganizing the terms of Eq.~\eqref{eq:CR3BP_J2000_transf} and Eq.~\eqref{eq:CR3BP_J2000_transf_v}, though this exercise is here omitted.

\subsection{Time correspondence}
\label{subsec:isochronous_correspondence}

In addition to being acquainted with the coordinate transformations discussed above, understanding the time dilation/contraction effects between the CR3BP and HFEM plays a major role in ensuring an adequate initial guess for the transition algorithms. These effects arise because the dynamics of the CR3BP are established in dimensionless form considering a fixed characteristic time that brings the angular velocity of the primaries to a constant unitary value, made possible only by the assumption of constant distance between the Earth and Moon. However, under the HFEM, the distance between the Earth and Moon is no longer constant. In this regard, directly translating time-related quantities from the CR3BP to HFEM will lead to a \textit{coherence} problem, given that this would break the unitary angular velocity assumption of the CR3BP, as pointed out in \cite{Park2024AssessmentDynModels}. While this is not necessarily impeding from the viewpoint of the trajectory transition process, it is evidently relevant to correct it for the purposes of benefiting numerical convergence. To this end, an extra step to translate the non-uniform correspondence between time in the CR3BP, which we may denote $\tilde{t}$, and time within the HFEM, $t$, is pursued to force the Earth--Moon line to rotate at unitary angular velocity, $n=1$, when viewed in the dimensionless time-frame. This is done by locally interpreting the orbits of the primaries as circles, such that the instantaneous variation of the Earth--Moon line angle, $\theta$, is obtained from the (common) expression for a circular orbit, i.e.
\begin{equation}
	\label{eq:time_corr}
	n = \frac{d\theta}{d\tilde{t}}= 1 \Longleftrightarrow \frac{d\theta}{dt} \frac{dt}{d\tilde{t}} = 1 \Longleftrightarrow \frac{d\tilde{t}}{dt} = \frac{d\theta}{dt}\Longleftrightarrow \frac{d\tilde{t}}{dt} = \sqrt{\frac{\mu_E+\mu_M}{d_\text{HFEM}^3}},
\end{equation}
Note that Eq.~\eqref{eq:time_corr} assumes that $\tilde{t}$ and $t$ are scaled by the same factor. This transformation ensures local time coherence between the two models' representation of the synodic reference frame.

Ultimately, Eq.~\eqref{eq:time_corr} may be integrated to more adequately translate the passing of time under the CR3BP to the HFEM. As such, a time instant $\tilde{t}_i$ under the CR3BP, relative to an established reference, is better represented in the HFEM by the time
\begin{equation*}
	t_i = \int_{0}^{\tilde{t}_i} \sqrt{\frac{d_\text{HFEM}^3}{\mu_E+\mu_M}} d\tilde{t}
\end{equation*}
relative to that same reference\footnote{Note that $d_\text{HFEM}$ is retrieved through ephemeris data resorting to dimensional time and it is hence necessary to scale $\tilde{t}$ by the CR3BP characteristic time and add the reference epoch in order to retrieve the corresponding instantaneous distance between the Earth and Moon.}.
Equivalently, any quantity relating to time that follows the typical CR3BP scaling ought to be adjusted according to Eq.~\eqref{eq:time_corr} before being passed to the HFEM, in order to translate the time dilation/contraction effects between the two dynamical models. To this end, the velocity transformation in Eq.~\eqref{eq:CR3BP_J2000_transf_v}, which considers the spacecraft velocity from the CR3BP as a derivative of position with respect to $\tilde{t}$, i.e. $\prescript{C}{}{\mathbf{v}}=d\prescript{C}{}{\mathbf{r}}/d\tilde{t}$, needs to be rewritten so as to adjust this contribution to be a derivative with respect to $t$ instead, i.e.
\begin{equation}
	\label{eq:CR3BP_J2000_transf_v_corr}
	\prescript{J}{}{\mathbf{v}} = \prescript{J}{}{\mathbf{v}}_b + \left(\frac{\dot{d}_\text{HFEM}}{d_\text{CR3BP}}\mathbf{R} + \frac{d_\text{HFEM}}{d_\text{CR3BP}}\dot{\mathbf{R}}\right)\prescript{C}{}{\mathbf{r}} + \frac{d_\text{HFEM}}{d_\text{CR3BP}}\mathbf{R}\frac{d\tilde{t}}{dt}\prescript{C}{}{\mathbf{v}}.
\end{equation}

\section{Trajectory design}
\label{sec:trajectory_design}

In literature, two main approaches to the problem of transitioning a trajectory from the CR3BP to the HFEM can typically be found. The simpler, so-called single-shooting methods, consider as an initial guess a single point in phase space that is iteratively adjusted until the propagation of the HFEM dynamics respects a set of objectives or operational constraints after a specified period of time, which may itself also be a variable of the problem. While conceptually simple, execution of single-shooting methods typically falters at retrieving accurate trajectories over large windows of time due to growing errors in the linear approximation of the highly nonlinear and chaotic HFEM dynamics \cite{spreen2021phd,soto2023thesis} --- a necessary step for the computation of the problem's Jacobian, used in most algorithms. Furthermore, the absence of control points along the integrated trajectory path may render the converged solutions undesirable or unfeasible to follow. In contrast, multiple-shooting methods divide a trajectory onto various segments that are fed to the transition algorithm, with the intent of alleviating the previous concerns. The spacecraft conditions at the beginning of each segment, so-called patch points, constitute the variables of the problem and should be adjusted not only until the set of objectives or operational constraints are met, but also to guarantee that there is continuity between the integration of successive segments. While these methods present a clear increase in terms of complexity, the benefits in numerical robustness and flexibility render multiple-shooting much better suited for the challenges of HFEM transition, reason as to why they were pursued in this work.

\subsection{Multiple-shooting}
\label{subsec:multiple_shooting}
\begin{figure}[h]
	\centering
	\includegraphics[width=.8\textwidth]{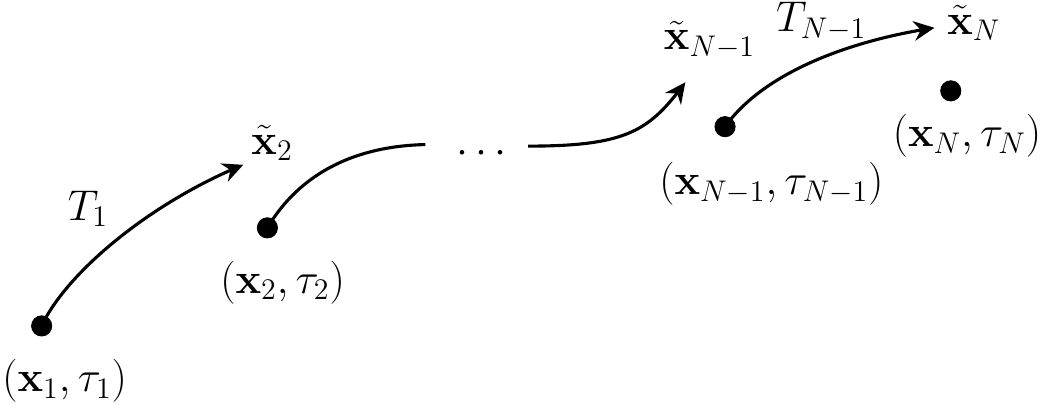}
	\caption{Multiple shooting example of a CR3BP initial guess under the HFEM.}
	\label{fig:multiple_shooting}
\end{figure}

As previously stated, trajectories in the CR3BP, deemed straightforward to retrieve, will serve as the initial guesses for the generation of equivalent solutions in the HFEM, through multiple-shooting. Naturally, under this more realistic dynamical model, the integration of the various segments that make up the CR3BP trajectory will not yield a continuous solution, as exemplified in Fig.~\ref{fig:multiple_shooting}. 
We denote by $\mathbf{x}_i^T = \begin{bmatrix}
	\mathbf{r}_i^T & \mathbf{v}_i^T
\end{bmatrix} \in \mathbb{R}^6, ~i = 1,\dots,N$, the states of the satellite at the $i$-th trajectory patch point, $\tau_i$ the epoch at which they occur (absolute time with respect to a fixed reference), and $T_i$ the time interval to the next patch point. Furthermore, $\tilde{\mathbf{x}}_{i+1}$ denotes the result of integrating the HFEM dynamics from the state $\mathbf{x}_i$, at epoch $\tau_i$, over the time interval $T_i$, i.e.
\begin{equation*}
	\tilde{\mathbf{x}}_{i+1} = \mathbf{x}_i +  \int_{\tau_i}^{\tau_i+T_i} \mathbf{f}(\mathbf{x}(t)) dt, \quad \text{with} \quad  \mathbf{x}(\tau_i) = \mathbf{x}_i,
\end{equation*}
where $\mathbf{f}(\mathbf{x})$ corresponds to the HFEM EoM, from Eq.~\eqref{eq:HFEM_eom}. Note that, since the dynamics are being propagated forward from each patch point, the problem is posed as forward multiple-shooting \cite{park2025CharacterizationL2Analogs,Spreen2023BaselineNRHO}. Alternatives such as forward-backward multiple-shooting \cite{williams2017NRHO, Diane2017stationkeeping} may provide benefits in convergence once the transition algorithms are established, but they typically represent an increase to the computational burden of the application and were hence discarded for the analysis presented in this work.

Given that the trajectory under the HFEM will initially be discontinuous, it is necessary to impose a \textit{continuity constraint} in order to ensure that the transition algorithms evolve towards a feasible solution. The nature of this continuity constraint depends on the multiple-shooting approach pursued, with two main alternatives found in literature. On the one hand, a fixed-time approach lets only the states of the spacecraft at each patch point vary, fixing the corresponding epochs and time intervals. Denoting the design variables by $\mathbf{X}$ and the continuity constraint by the function $\mathbf{F}(\mathbf{X})$, a fixed-time formulation therefore considers
\begin{equation}
	\label{eq:op_prob_fixed}
	\mathbf{X} = \begin{bmatrix}
		\mathbf{x}_1 \\
		\vdots\\
		\mathbf{x}_N
	\end{bmatrix}
	\qquad \text{and} \qquad
	\mathbf{F}(\mathbf{X}) = \begin{bmatrix}
		\tilde{\mathbf{x}}_2 - \mathbf{x}_2 \\
		\vdots \\
		\tilde{\mathbf{x}}_N - \mathbf{x}_N
	\end{bmatrix}.
\end{equation}
In contrast, in a variable-time approach the epochs of each patch point and/or the time intervals between them accompany the spacecraft states as design variables of the transition algorithm \cite{park2025CharacterizationL2Analogs}. In theory, this alternative is more versatile and better suited for the time-varying dynamics of the HFEM. However, the inclusion of additional variables in $\mathbf{X}$ and constraints in $\mathbf{F}(\mathbf{X})$ typically constitutes a meaningful increase in terms of complexity --- both at the theoretical and implementation level --- given that it is usually necessary to compute the Jacobian matrix to proceed numerically, as will be seen shortly. In fact, while variable-time formulations introduce additional degrees of freedom that may facilitate meeting additional operational constraints or objectives, it is not necessarily the case that numerical convergence is improved. In this sense, the benefits in implementation tractability and computational burden provided by fixed-time formulations of the shooting problem often make them an optimal choice in a first analysis. Moreover, preliminary tests have shown that the absence of time-related quantities in the design vector is often inconsequential if efforts are put forth towards assuring an adequate time correspondence between the CR3BP and HFEM, as covered in Section~\ref{subsec:isochronous_correspondence}. To this end, only a fixed-time formulation has been considered for the analysis presented in this work.

\subsection{Additional constraints}
\label{sec:add_constraints}
In addition to ensuring that there is continuity in the converged solution, it is also desirable to have some degree of control over its overall shape. Considering that HFEM trajectories tend to exist near equivalent CR3BP structures, which most transition schemes take advantage of, coercing the algorithms to remain in the vicinity of the initial CR3BP guess may be particularly relevant for boosting convergence and robustness. Moreover, since the CR3BP solutions are typically chosen carefully for optimal positioning with respect to the celestial bodies or communications clearance, for example, this concern is also relevant from an operational perspective. In this sense, the inclusion of a constraint that dictates proximity to the original CR3BP solution being transitioned into the HFEM is expected to be a valuable addition in most cases.

For the sake of generality, we consider the case where the entire vector of design variables is to lie within some measure of proximity to the initial guess, which therefore constitute \textit{desired conditions}, denoted by $\mathbf{X}_d$. To this end, we define the constraint
\begin{equation}
	\label{eq:prox_constraint}
	\norm{\mathbf{X}-\mathbf{X}_d}^2_\mathbf{W} \leq \rho,
\end{equation}
where the constant $\rho>0$ tunes how close the converged solution should be to the desired conditions and $\mathbf{W}\in\mathbb{S}_{++}$ is a symmetric positive definite matrix, allowing one to select which entries of the design variables are meaningful for the proximity objective --- for example, we might be concerned only with the position entries of $\mathbf{X}$ or the initial/final states of the spacecraft along its trajectory. For the sake of intuition, we consider $\mathbf{W}$ to be a diagonal matrix.

In addition to the proximity constraint, other operational objectives depending on specific requirements of each particular mission may be pondered and included either as hard constraints at the level of $\mathbf{F}(\mathbf{X})$, or via the use of inequalities, similarly to Eq.~\eqref{eq:prox_constraint}. To keep this analysis as general as possible, such concerns are not going to be considered at this point.

\section{Optimization problem formulation}
\label{sec:op_probs}
Transitioning trajectories from the CR3BP to the HFEM via multiple-shooting solely under the continuity constraint corresponds to finding a solution of the nonlinear equation $\mathbf{F}(\mathbf{X})=\mathbf{0}$, which encodes the propagation of the HFEM dynamics over each trajectory segment, as detailed in Eq.~\eqref{eq:op_prob_fixed}. Since this problem is under-constrained, we expect multiple solutions to be available. Furthermore, given its nonlinear nature, there is no hope for retrieving a general solution in a trivial manner. Hence, it is logical to pursue an optimization approach to tackle the problem of meeting the continuity constraint, which should also facilitate the inclusion of other relevant objectives, such as the proximity goal detailed in Section~\ref{sec:add_constraints}.

In literature, the optimization problems are typically constructed considering the continuity objective directly at the level of the constraints, i.e. through a formulation in the likeness of
\begin{equation}
	\label{eq:feas_prob}
	\begin{aligned}
		\min_\mathbf{X}& \quad g(\mathbf{X})\\
		\mathrm{s.t.}& \quad \mathbf{F}(\mathbf{X}) = \mathbf{0},
	\end{aligned}
\end{equation}
where $g(\mathbf{X})$ is a scalar function that depends on the particular implementation followed. Given that $\mathbf{F}(\mathbf{X})$ is nonlinear, the resulting optimization problem is non-convex. A common approach to address this concern is to approximate the original optimization problem via a sequence of convex counterparts. In this sense, while the formulation in Eq.~\eqref{eq:feas_prob} may represent the most literal interpretation of the continuity goal, we will later see that carrying out such an approximation at the level of a hard constraint introduces relevant concerns from a theoretical standpoint.

In contrast, we consider a relaxed version of the optimization problem that takes the continuity goal to the objective function of the optimization problem, through a least-squares viewpoint. If the proximity constraint is also included, an approach of this nature takes the form of
\begin{equation}
	\label{eq:prox_const_prob}
	\begin{aligned}
		\min_\mathbf{X}& \quad \norm{\mathbf{F}(\mathbf{X})}^2\\
		\mathrm{s.t.}& \quad \norm{\mathbf{X}-\mathbf{X}_d}_\mathbf{W}^2 \leq \rho.
	\end{aligned}
\end{equation}
Alternatively, we may establish the optimization problem in a fully unconstrained manner as
\begin{equation}
	\label{eq:unconst_prob}
	\begin{aligned}
		\min_\mathbf{X}& \quad \norm{\mathbf{F}(\mathbf{X})}^2 + \norm{\mathbf{X}-\mathbf{X}_d}_\mathbf{Q}^2,
	\end{aligned}
\end{equation}
where $\mathbf{Q}\in\mathbb{S}_{++}$ is appropriately chosen to guarantee the original proximity constraint in Eq.~\eqref{eq:prox_constraint} --- a choice which might, however, be nontrivial to make and will shortly warrant further discussion. Nonetheless, the fully unconstrained version of the optimization problem will prove ideal for the determination of an analytical update step once convexification is pursued. Moreover, we will see that it allows for the construction of an algorithm that ensures $\mathbf{F}(\mathbf{X})$ is non-increasing between successive iterations, bestowing the trajectory design process with a measure of robustness and stability.

\section{Convexification}
\label{sec:convex}
To have any chance at generally solving the optimization problems detailed in Section~\ref{sec:op_probs}, it is necessary to address the non-convex nature of the continuity objective encoded through the nonlinear function $\mathbf{F}(\mathbf{X})$, a key element of all formulations presented insofar. In this work, the original optimization problems are tackled through a sequence of convex approximations, each established about a nominal condition corresponding to the value of the design variables at the previous iteration. Denoting the design vector at the $k$-th iteration by $\mathbf{X}_k$, we take $\mathbf{X} = \mathbf{X}_k + \Delta \mathbf{X}_k$, such that, for a sufficiently small $\Delta \mathbf{X}_k$, the nonlinear function may be approximated through
\begin{equation}
	\label{eq:F_approx}
	\mathbf{F}(\mathbf{X}) \approx \mathbf{F}_k + \mathbf{J}_k\Delta \mathbf{X}_k,
\end{equation}
where $\mathbf{F}_k:=\mathbf{F}(\mathbf{X}_k)$ is the function evaluated at $\mathbf{X}_k$ and $\mathbf{J}_k:=\mathbf{J}_\mathbf{F}(\mathbf{X}_k)$ is the Jacobian of $\mathbf{F}(\mathbf{X})$ evaluated at the same point. As further detailed in \cite{spreen2021phd}, the Jacobian for the fixed-time formulation of the transition problem under consideration (i.e., the one in Eq.~\eqref{eq:op_prob_fixed}) is given by
\begin{equation}
	\label{eq:jac}
	\mathbf{J}_\mathbf{F}(\mathbf{X}) = \begin{bmatrix}
		\frac{\partial \tilde{\mathbf{x}}_2}{\partial \mathbf{x}_1} & - \frac{\partial \mathbf{x}_2}{\partial \mathbf{x}_2} & \\[1.5ex]
		& \ddots & \ddots &\\[2ex]
		& & \frac{\partial \tilde{\mathbf{x}}_N}{\partial \mathbf{x}_{N-1}} &  -\frac{\partial \mathbf{x}_N}{\partial \mathbf{x}_{N}}
	\end{bmatrix} = 
	\begin{bmatrix}
		\mathbf{\Phi}_1 & -\mathbf{I_6} & \\[1ex]
		& \ddots & \ddots &\\[1.5ex]
		& & \mathbf{\Phi}_{N-1} & -\mathbf{I_6}
	\end{bmatrix},
\end{equation}
where $\mathbf{\Phi}_j=\mathbf{\Phi}(\tau_j+T_j,\tau_{j})$ is the state-transition matrix (STM) over the $j$-th trajectory segment, and empty entries are zero. The STM is retrieved by integrating the HFEM dynamics in Eq.~\eqref{eq:HFEM_eom}, after linearization. This standard, \textit{variational} approach is preferred over a numerical alternative (e.g., via finite differences) due to the inherent sensitivity of the dynamics at play. Despite the fact that $\mathbf{J}_k$ is exact, the linear approximation behind Eq.~\eqref{eq:F_approx} motivates the use of multiple-shooting approaches over single-shooting alternatives, since subdivision of long trajectories into shorter segments improves the validity of this approximation. Note that $\mathbf{J}_k$ is always full row-rank.

The solution of each convex optimization problem, under the approximation in Eq.~\eqref{eq:F_approx}, yields the optimal step $\Delta \mathbf{X}_k^*$ that updates the set of design variables from the previous iteration, i.e. $\mathbf{X}_{k+1} = \mathbf{X}_k + \Delta \mathbf{X}_k^*$.
In this sense, note that in a least-squares formulation where $\norm{\mathbf{F}(\mathbf{X})}^2$ is the sole contribution to the objective function of the optimization problem, such as in Eq.~\eqref{eq:prox_const_prob}, the optimal step guarantees, by definition, that its approximate value is non-increasing, i.e.
\begin{equation*}
	\norm{\mathbf{F}_{k} + \mathbf{J}_k\Delta \mathbf{X}_k^*}^2  \leq \norm{\mathbf{F}_k}^2.
\end{equation*}
In addition, if the linear approximation is accurate, or $\Delta \mathbf{X}_k^*$ is sufficiently small, this effect is reflected in the original function, i.e.
\begin{equation*}
	\norm{\mathbf{F}\left(\mathbf{X}_k + \Delta \mathbf{X}_k^*\right) }^2\approx \norm{\mathbf{F}_{k} + \mathbf{J}_k\Delta \mathbf{X}_k^*}^2 \leq \norm{\mathbf{F}(\mathbf{X}_k)}^2.
\end{equation*}
This makes the relaxed versions of the optimization problem particularly reliable, in contrast with formulations where $\mathbf{F}(\mathbf{X})=\mathbf{0}$ is taken at the constraint-level, much like in Eq.~\eqref{eq:feas_prob}. In fact, one might argue that formulations of this nature are badly posed from a theoretical perspective, given how meeting the approximate constraint $\mathbf{F}_k + \mathbf{J}_k \Delta \mathbf{X}_k = \mathbf{0}$ will generally not yield $\mathbf{F}(\mathbf{X}_k+\Delta \mathbf{X}_k) = \mathbf{0}$.

\subsection{Minimum-norm update}
\label{subsec:MN}
Despite the previous concern, most of the literature opts for approaches based on formulations in the likeness of Eq.~\eqref{eq:feas_prob} for their ease of implementation or versatility. To potentially mitigate the shortcomings discussed, it is typical to minimize the norm of the update step at each iteration, i.e. solving a sequence of optimization problems in the form of
\begin{equation}
	\label{eq:min_norm_prob}
	\begin{aligned}
		\min_{\Delta \mathbf{X}_k}& \quad \norm{\Delta \mathbf{X}_k}^2\\
		\mathrm{s.t.}& \quad \mathbf{F}_k + \mathbf{J}_k\Delta \mathbf{X}_k = \mathbf{0},
	\end{aligned}
\end{equation}
which correspond to quadratic programs (QPs). Being a QP with only a linear equality constraint means that the problem in Eq.~\eqref{eq:min_norm_prob} may be solved directly to obtain an analytical expression for the optimal update step. To do so, the auxiliary variable $\mathbf{\lambda}$ is introduced to encode the equality constraint and write the Lagrangian
$\mathcal{L} = \norm{\Delta \mathbf{X}_k}^2 + \boldsymbol{\lambda}^T(\mathbf{F}_k + \mathbf{J}_k\Delta \mathbf{X}_k)$.
Then, taking the gradients with respect to $\Delta \mathbf{X}_k$ and $\boldsymbol{\lambda}$ to zero yields the so-called MN update step equation, i.e.
\begin{equation*}
	\begin{cases}
		\nabla_{\Delta \mathbf{X}_k} \mathcal{L} = 2 \Delta \mathbf{X}_k + \mathbf{J}_k^T \boldsymbol{\lambda} = 0 \implies \Delta \mathbf{X}_k = - \frac{1}{2}\mathbf{J}_k^T\boldsymbol{\lambda} & \\
		\nabla_{\boldsymbol{\lambda}} \mathcal{L} = \mathbf{J}_k\Delta \mathbf{X}_k + \mathbf{F}_k = 0 &
	\end{cases}
\end{equation*}
leads to
\begin{equation}
	\label{eq:MN_step}
	\boldsymbol{\lambda} = 2 \left(\mathbf{J}_k \mathbf{J}_k^T\right)^{-1}\mathbf{F}_k \implies \Delta \mathbf{X}^*_k = - \mathbf{J}_k^T\left(\mathbf{J}_k \mathbf{J}_k^T\right)^{-1} \mathbf{F}_k.
\end{equation}
Despite the intent of improving the quality of the linear approximation through the minimization of the update step, Eq.~\eqref{eq:MN_step} does not guarantee the decrease of $\norm{\mathbf{F}(\mathbf{X})}$, and the iterative procedure is susceptible to divergence. Some authors also raise concerns on its efficacy at tackling large scale problems \cite{williams2017NRHO}. For these reasons, it is often necessary to limit the magnitude of the update step, such that its norm does not surpass a given threshold, as suggested in \cite{Spreen2023BaselineNRHO}. In other words, the update step is modified according to
\begin{equation}
	\label{eq:MN_attenuation}
	\Delta \mathbf{X}^\star_k = \Delta \mathbf{X}^*_k \frac{\min\left(\gamma, \norm{\Delta \mathbf{X}^*_k}\right)}{\norm{\Delta \mathbf{X}^*_k}},
\end{equation}
where $\gamma>0$ is a constant selected manually, depending on the application. However, this technique acts as a way of minimizing the risks of divergence rather than fully addressing them. Moreover, as the magnitude threshold is selected manually, the value that provides the best trade-off between robustness and convergence speed is not necessarily obvious or general, as we will later demonstrate.

Furthermore, note that the inclusion of the proximity objective in this formulation of the optimization problem may be problematic. If it were to figure the objective function in a relaxed manner, much like Eq.~\eqref{eq:unconst_prob}, this contribution could potentially overthrow the minimization of the update step norm, putting into question the numerical stability of the algorithm. On the contrary, if it is brought to the constraints of the problem in the form of an inequality, as in Eq.~\eqref{eq:prox_const_prob}, it would make it cumbersome to compute an analytical update equation, rendering this approach significantly more computationally demanding. For these reasons, proximity objectives in MN techniques from literature are typically posed as additional equality constraints, by augmenting $\mathbf{F}(\mathbf{X})$ \cite{pavlak2012Strategy,pavlak2010thesis}. While simple objectives may be met in this manner, it is also the case that the problem easily becomes over-constrained if multiple goals are mandated, noticing how $\mathbf{X}$ and $\mathbf{F}(\mathbf{X})$ differ in size only by 6 entries. For example, feasibility issues quickly arise if it is necessary to impose various constraints in terms of three-dimensional position at different specified epochs, to reflect operational objectives. In such cases, a regular fixed-time MN technique is unable to proceed, and a reformulation of the shooting problem considering variable-time segments is typically needed to introduce additional degrees of freedom, as discussed at the end of Section~\ref{subsec:multiple_shooting}. This, however, represents a step up in implementation complexity. Moreover, the additional room for constraints may be quickly populated if specific patch points are constrained in epoch, e.g. to avoid eclipsing. For these reasons, MN-based approaches may not provide sufficient control over the converged solution's shape, especially under a fixed-time formulation --- a concern to be evaluated later, when numerical trials are performed in Section~\ref{sec:results}.

Despite of the theoretical shortcomings discussed above, the MN update equation is, in many cases, very capable of transitioning solutions from the CR3BP to the HFEM. As such, it is at the base of many algorithms found in the literature for its benefits in ease of implementation and fast convergence speeds \cite{park2025CharacterizationL2Analogs,Spreen2020NRHO,Park2024AssessmentDynModels}. Due to its widespread influence in literature, the MN update step was thus selected as a baseline for the evaluation of the alternative developed in the remainder of this work. Its implementation is laid out in Algorithm~\ref{alg:MN}. Convergence is assumed to be reached once $\norm{\mathbf{F}_k}\leq \sigma$, with $\sigma>0$ a desired tolerance. To prevent the algorithm from running indefinitely, the number of iterations is limited to $k_{\text{max}}>0$. In practice, it may be necessary to employ additional stopping criteria to terminate execution in the case of divergence.

\begin{algorithm}
	\caption{MN update scheme with step attenuation.}
	\label{alg:MN}
	\begin{algorithmic}[1]
		\Require $\mathbf{X}_0$, $k_{\text{max}}$, $\sigma$, $\gamma$
		\Ensure Final solution, $\bar{\mathbf{X}}$
		\State Compute $\mathbf{F}_0$, $\mathbf{J}_0$
		\State Initialize $k=0$
		\Repeat
		\State Compute $\Delta \mathbf{X}_k^\star$ via Eqs.~\eqref{eq:MN_step} and \eqref{eq:MN_attenuation}
		\State Set $\mathbf{X}_{k+1} = \mathbf{X}_k+\Delta \mathbf{X}_k^\star$
		\State Compute $\mathbf{F}_{k+1}$, $\mathbf{J}_{k+1}$
		\State Increment $k\gets k+1$
		\Until{$\norm{\mathbf{F}_k}<\sigma$ \textbf{or} $k>k_{\text{max}}$}
		\State Retrieve $\bar{\mathbf{X}}=\mathbf{X}_k$
	\end{algorithmic}
\end{algorithm}

\subsection{Quadratically constrained quadratic programs}
Given the previous discussion on the benefits of the relaxed formulations of the multiple-shooting problem, this work investigates the formulations in Eq.~\eqref{eq:prox_const_prob} and Eq.~\eqref{eq:unconst_prob} as an alternative to Eq.~\eqref{eq:feas_prob}, explored in the previous section. In this sense, the optimization problem in Eq.~\eqref{eq:prox_const_prob} is approximated by a sequence of convex counterparts in the form of
\begin{equation}
	\label{eq:prox_const_LS}
	\begin{aligned}
		\min_{\Delta \mathbf{X}_k}& \quad \Delta\mathbf{X}_k^T\mathbf{J}_k^T\mathbf{J}_k\Delta\mathbf{X}_k + 2\Delta \mathbf{X}_k^T\mathbf{J}_k^T\mathbf{F}_k\\
		\mathrm{s.t.}& \quad \norm{\Delta \mathbf{X}_k + \mathbf{E}_k}_\mathbf{W}^2 \leq \rho,
	\end{aligned}
\end{equation}
after expanding and dropping constant terms of the objective function, where $\mathbf{E}_k:=\mathbf{X}_k-\mathbf{X}_d$ denotes the error between the $k$-th design vector and the desired conditions. The formulation in Eq.~\eqref{eq:prox_const_LS} corresponds to a quadratically constrained quadratic program (QCQP), a class of convex optimization problems for which, unfortunately, a closed form solution cannot be generally found. This means that more advanced techniques need to be employed in order to compute $\Delta \mathbf{X}_k^*$ at each iteration, which would increase the computational burden of the implementation. With this in mind, the QCQP approach was disregarded in this work, being left for future investigation.

\subsection{Levenberg--Marquardt algorithm}
\label{subsec:LM}
Similarly to the case covered in the previous section, the convexification of the formulation in Eq.~\eqref{eq:unconst_prob} leads to a sequence of optimization problems in the form of
\begin{equation}
	\label{eq:unconst_LS}
	\begin{aligned}
		\min_{\Delta \mathbf{X}_k} & \quad \Delta\mathbf{X}_k^T\left(\mathbf{J}_k^T\mathbf{J}_k + \mathbf{Q}\right)\Delta\mathbf{X}_k + 2\Delta \mathbf{X}_k^T\left(\mathbf{J}_k^T\mathbf{F}_k + \mathbf{Q}\mathbf{E}_k\right),
	\end{aligned}
\end{equation}
which correspond to QPs in $\Delta \mathbf{X}_k$, much like Eq.~\eqref{eq:min_norm_prob}. With this in mind, it is possible to obtain an analytical expression for the optimal step $\Delta \mathbf{X}_k^*$ that solves Eq.~\eqref{eq:unconst_LS} which, in this case, may be done directly by taking the gradient of the objective function to zero, leading to the well-known Newton--Gauss update step,
\begin{equation}
	\label{eq:NG_step}
	\begin{aligned}
		&\nabla\left[\Delta\mathbf{X}_k^T\left(\mathbf{J}_k^T\mathbf{J}_k + \mathbf{Q}\right)\Delta\mathbf{X}_k + 2\Delta \mathbf{X}_k^T\left(\mathbf{J}_k^T\mathbf{F}_k + \mathbf{Q}\mathbf{E}_k\right)\right] = 0\\  \Longleftrightarrow~&  2 \left(\mathbf{J}_k^T\mathbf{J}_k+\mathbf{Q}\right)\Delta \mathbf{X}_k^* + 2 \left(\mathbf{J}_k^T\mathbf{F}_k+\mathbf{Q}\mathbf{E}_k\right)= 0\\ \Longleftrightarrow~& \Delta \mathbf{X}_k^* = - \left(\mathbf{J}_k^T\mathbf{J}_k + \mathbf{Q}\right)^{-1}\left(\mathbf{J}_k^T\mathbf{F}_k + \mathbf{Q}\mathbf{E}_k\right).
	\end{aligned}
\end{equation}
Given the relaxed nature of the approach pursued, the proximity objectives are naturally included in the update step, without requiring a modification to the Jacobian matrix.
However, since $\mathbf{Q}$ is not necessarily full-rank and $\mathbf{J}_k$ is generally not square, there are no guarantees for the invertibility of $\mathbf{J}_k^T\mathbf{J}_k + \mathbf{Q}$. Hence, direct application of Eq.~\eqref{eq:NG_step} is typically not possible.

As a workaround that ultimately proves to be insightful, we consider the inclusion of a damping term in the objective function of the optimization problem, i.e. Eq.~\eqref{eq:unconst_LS} is rewritten as
\begin{equation*}
	\min_{\Delta \mathbf{X}_k} \quad \Delta\mathbf{X}_k^T\left(\mathbf{J}_k^T\mathbf{J}_k + \mathbf{Q}\right)\Delta\mathbf{X}_k + 2\Delta \mathbf{X}_k^T\left(\mathbf{J}_k^T\mathbf{F}_k + \mathbf{Q}\mathbf{E}_k\right) + \beta_k \Delta\mathbf{X}_k^T \Delta \mathbf{X}_k,
\end{equation*}
where $\beta_k>0$ is an iteration-dependent scalar that weights the square-norm of the step, $\Delta \mathbf{X}_k$. Once more, by taking the gradient of the objective function it is possible to arrive to an analytical update equation in the form of
\begin{equation}
	\label{eq:LM_step}
	\Delta \mathbf{X}_k^* = - \left(\mathbf{J}_k^T\mathbf{J}_k + \mathbf{Q} + \beta_k \mathbf{I}\right)^{-1}\left(\mathbf{J}_k^T\mathbf{F}_k + \mathbf{Q}\mathbf{E}_k\right),
\end{equation}
where $\mathbf{I}$ is the identity matrix with appropriate dimensions.
On the one hand, the additional term $\beta_k\mathbf{I}$ immediately solves the invertibility issue, ensuring that Eq.~\eqref{eq:LM_step} is always feasible. In addition, it provides direct control over the step size, which is evidently relevant given the linear approximation at the basis of the convexification. In this sense, the damping term $\beta_k$ can be increased to promote smaller update steps when the linear approximation is deemed poor, tending to a gradient-descent update, i.e.
\begin{equation*}
	\beta_k \to \infty \implies \Delta \mathbf{X}_k^* \approx -\frac{1}{\beta_k}\left(\mathbf{J}_k^T\mathbf{F}_k + \mathbf{Q}\mathbf{E}_k\right).
\end{equation*}
On the contrary, if the linear approximation is found to be adequate, $\beta_k$ may be decreased to improve convergence speed, tending to the Newton--Gauss update in Eq.~\eqref{eq:NG_step} when $\beta_k\to0$.

The standard LM algorithm corresponds to a systematic procedure to the definition of the damping term, based on the evolution of the original objective function, in this case from Eq.~\eqref{eq:unconst_prob}. In this work, we modify slightly the classical LM description to adjust $\beta_k$ based purely on the continuity residual, following two basic principles: \textit{(i)} if $\norm{\mathbf{F}_k}\leq \norm{\mathbf{F}_{k-1}}$, then the current step should be accepted and the damping term is updated for the next iteration following $\beta_{k+1}=\alpha \beta_k$, with $0<\alpha<1$; \textit{(ii)} if $\norm{\mathbf{F}_k}> \norm{\mathbf{F}_{k-1}}$, then the current step should be iteratively recalculated with a new damping term $\beta_k \gets \eta \beta_k$, with $\eta>1$, until the first condition is met. The values of $\alpha$ and $\eta$ are selected to promote convergence speed while avoiding unnecessary oscillations of the damping factor --- appropriate values for these parameters are discussed later, in Section~7, once specific application cases are targeted. By automatically controlling the update step, the algorithm requires no manual \textit{a posteriori} adjustments that a typical MN approach demands, such as the magnitude control detailed in Eq.~\eqref{eq:MN_attenuation}. Moreover, employing  this modified LM scheme evidently guarantees that $\norm{\mathbf{F}(\mathbf{X})}$ is non-increasing over successive iterations. This imparts a measure of robustness and stability to the algorithm that also distinguishes it from the MN update step in Eq.~\eqref{eq:MN_step}. Still, it cannot be argued at this point if the monotonous nature of the proposed scheme is a benefit or a disadvantage, as this particularity may be extremely useful in some scenarios but too stringent in others, where the algorithm may converge towards trajectories that identify local minima of $\norm{\mathbf{F}(\mathbf{X})}^2$ for which the continuity constraint is not met to the desired degree. To this end, numerical investigation of the proposed algorithm as an alternative to the MN update equation is warranted, being left for Section~\ref{sec:results}.

For the sake of simplicity, the implementation proposed in this work is termed simply ``LM algorithm'' from here on, and follows the listing in Algorithm~\ref{alg:LM}. The scheme is terminated once convergence is reached, according to the specified threshold on the continuity residual, or the maximum number of iterations is surpassed. Similarly to the MN alternative in Algorithm~\ref{alg:MN}, it may be necessary to impose additional stopping criteria. In this case, these should account for potential stalling at local minima, namely by limiting the execution of the inner LM loop. This is later discussed in Section~\ref{sec:algorithm}.
\begin{algorithm}
	\caption{LM algorithm with proximity objective.}
	\label{alg:LM}
	\begin{algorithmic}[1]
		\Require $\mathbf{X}_0$, $\mathbf{Q}$, $k_{\text{max}}$, $\sigma$, $\beta_0$, $\alpha$, $\eta$
		\Ensure Final solution, $\bar{\mathbf{X}}$
		\State Compute $\mathbf{F}_0$, $\mathbf{J}_0$
		\State Initialize $k=0$
		\Repeat \Comment{Outer loop}
		\Repeat \Comment{Inner LM loop}
		\State Compute $\Delta \mathbf{X}_k^*$ via Eq.~\eqref{eq:LM_step}
		\State Try $\mathbf{X}_{k+1} = \mathbf{X}_k+\Delta \mathbf{X}_k^*$
		\State Compute $\mathbf{F}_{k+1}$, $\mathbf{J}_{k+1}$
		\If{$\norm{\mathbf{F}_{k+1}}>\norm{\mathbf{F}_{k}}$}
		\State Update $\beta_k\gets \eta \beta_k$ \Comment{Refuse step, increase damping}
		\EndIf
		\Until{$\norm{\mathbf{F}_{k+1}}\leq\norm{\mathbf{F}_{k}}$}
		\State Update $\beta_{k+1} = \alpha \beta_k$ \Comment{Accept step, decrease damping}
		\State Increment $k\gets k+1$
		\Until{$\norm{\mathbf{F}_{k}}<\sigma$ \textbf{or} $k>k_\text{max}$}
		\State Retrieve $\bar{\mathbf{X}}=\mathbf{X}_k$
	\end{algorithmic}
\end{algorithm}

\subsection{Extensions}
\label{subsec:extensions}
Recall from Eq.~\eqref{eq:unconst_prob} that proximity to the CR3BP guess is weighted under the diagonal matrix $\mathbf{Q}$, which is reflected in the proposed LM update step in Eq.~\eqref{eq:LM_step}. While it may be evident to decide which entries of $\mathbf{Q}$ should be non-zero to meet specific goals, the determination of the optimal values for these entries is non-trivial, since this term may significantly shape the objective function being minimized. To this end, the use of an adaptive procedure to the definition of the entries of $\mathbf{Q}$ is relevant. Alternatively, since the minimization of the objective function is only concerned with the relative weight between the proximity and continuity goals, we instead propose to shape the contribution $\norm{\mathbf{F}(\mathbf{X})}^2$. In this sense, we restate the original least-squares formulation from Eq.~\eqref{eq:unconst_prob} as
\begin{equation}
	\label{eq:G_unconstrained_least_squares}
	\min_{\mathbf{X}} \quad G\left(\norm{\mathbf{F}(\mathbf{X})}^2\right) + \norm{\mathbf{X}-\mathbf{X}_{d}}_\mathbf{Q}^2,
\end{equation}
where $G(u):\mathbb{R}\to\mathbb{R}$ is a scalar function. Ideally, $G(u)$ should decrease significantly as $u$ is brought to zero, in order to further incentivize meeting the continuity objective. In this work, we consider the concave function
\begin{equation*}
	G(u) = \log(\sqrt{u+\delta} + \epsilon),
\end{equation*}
where $0<\delta,\epsilon\ll 1$ are tuning parameters, whose definition is, once more, not entirely evident \textit{a priori}. If too large, decreasing $u$ yields no benefits, and if too small, numerical instability might be introduced for incentivizing $u$ to decrease too hastily. Hence, an adaptive formulation to the definition of the parameters is also pursued, namely taking $\delta$ and $\epsilon$ to be $n$ orders of magnitude below $\norm{\mathbf{F}}^2$ and  $\norm{\mathbf{F}}$, respectively --- a process known as \textit{annealing}. Preliminary testing has shown that $n=2$ yields adequate results.

Since $G(u)$ is concave, the optimization problem in Eq.~\eqref{eq:G_unconstrained_least_squares} cannot be solved directly in a general manner, meaning we once more resort to solving a sequence of convex approximations, linearized about the converged solution of the previous optimization problem. After dropping constant terms, each optimization problem may thus be written as
\begin{equation}
	\label{eq:adap_unconst_LS}
	\min_{\mathbf{X}} \quad \xi\norm{\mathbf{F}(\mathbf{X})}^2 + \norm{\mathbf{X}-\mathbf{X}_{d}}_\mathbf{Q}^2,
\end{equation}
where the adaptive weight, $\xi$, is given by
\begin{equation}
	\label{eq:kappa}
	\xi := \dot{G}\left(\norm{\mathbf{F}(\bar{\mathbf{X}})}^2\right) = \frac{1}{2\sqrt{\norm{\mathbf{F}(\bar{\mathbf{X}})}^2+\delta}\left(\sqrt{\norm{\mathbf{F}(\bar{\mathbf{X}})}^2+\delta}+\epsilon\right)},
\end{equation}
with $\bar{\mathbf{X}}$ denoting the converged solution of the previous optimization problem. Each optimization problem is non-convex due to $\mathbf{F}(\mathbf{X})$, and is thus solved internally via the proposed LM algorithm. When doing so, it is important to remember to propagate the weight $\xi$ to the update equation used in each $k$-th inner iteration, which becomes
\begin{equation}
	\label{eq:LM_step_adap_weight}
	\Delta \mathbf{X}_k^* = -\left(\xi \mathbf{J}_k^T\mathbf{J}_k + \mathbf{Q} + \beta_k \mathbf{I}\right)^{-1}\left(\xi \mathbf{J}_k^T\mathbf{F}_k + \mathbf{Q} \mathbf{E}_k\right).
\end{equation}
The same algorithmic logic for the definition of damping the term, $\beta_k$, presented at the end of Section~\ref{subsec:LM}, should be pursued.

Note how the formulation in Eq.~\eqref{eq:adap_unconst_LS} follows closely that of the original unconstrained least-squares problem in Eq.~\eqref{eq:unconst_LS}, except for the additional adaptive weight that automatically dictates whether the proximity or continuity goal dominates. In this work, we initially take $\xi=1$ to weight these two contributions equally and, if convergence is subpar, employ the annealing process detailed above until the continuity residual meets a desired threshold.
Ultimately, the gradual process of increasing the preponderance of the continuity objective is expected to result in a final converged solution that meets the proximity objective as best as possible while bringing $\norm{\mathbf{F}}$ within the required threshold.

\section{Final algorithm}
\label{sec:algorithm}

\begin{figure}[h]
	\centering
	\begin{minipage}{.48\linewidth}
		\centering
		\includegraphics[height=10cm]{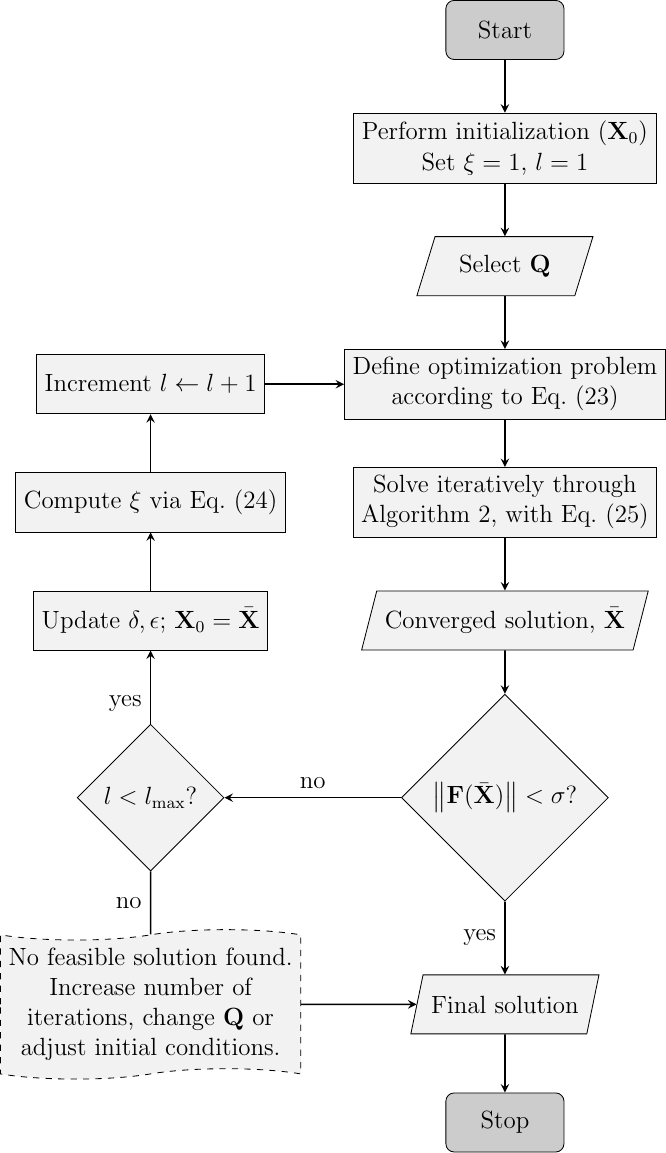}
		\caption{LM implementation considering adaptive weighting.}
		\label{fig:flowchart_LM}
	\end{minipage}
	\hfill
	\begin{minipage}{.48\linewidth}
		\centering
		\includegraphics[height=10cm]{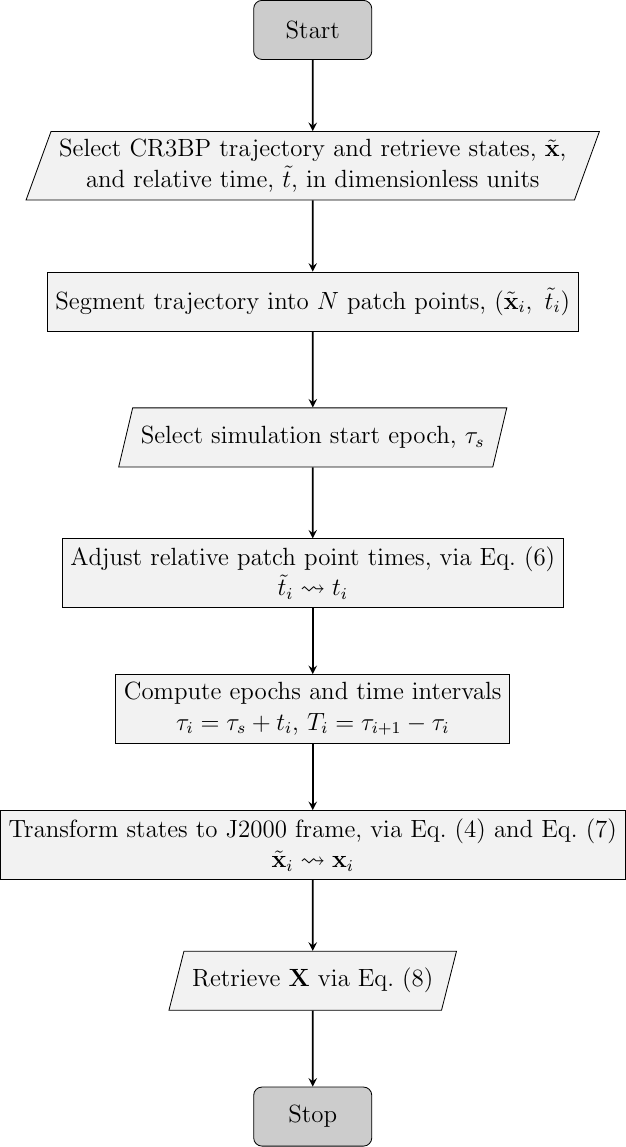}
		\caption{Construction of the initial guess from a CR3BP solution.}
		\label{fig:flowchart_init}
	\end{minipage}
\end{figure}

The overall LM workflow is presented schematically in Fig.~\ref{fig:flowchart_LM}, considering the possibility of employing adaptive weighting in the formulation of the optimization problem --- which may be discarded by setting $l_\text{max}=1$. An initial factor $\xi=1$ is selected, in order to attribute equal weights to the continuity and proximity objectives. As previously detailed in Algorithm~\ref{alg:LM}, the LM scheme used to solve each optimization problem is stopped once the maximum number of iterations has been reached or $\norm{\mathbf{F}}$ falls within a desired threshold, $\sigma$, in which case a final converged solution has been achieved. In practice, however, it may also be necessary to terminate execution if the relative variation of $\norm{\mathbf{F}}$ between successive iterations is small, since this points towards stalling at a local minimum. This is especially relevant if adaptive weighting is being employed, as it suggests that $\xi$ may need to be adjusted to potentially achieve a solution that brings $\norm{\mathbf{F}}$ below $\sigma$. If this is not possible after solving $l\geq l_\text{max}$ optimization problems, then the local minimum is assumed to be persistent and the algorithm is deemed unable to bring the continuity residual to a sufficiently satisfactory value. This may occur due to an insufficient number of iterations, inadequate definition of $\mathbf{Q}$ or the thresholds for the stopping criteria, or a poor initial guess --- the latter of which may be indicative of a transition-challenging region.

The process of constructing an initial guess from a CR3BP solution is presented in Fig.~\ref{fig:flowchart_init}, and is irrespective of the transition algorithm employed. This procedure follows the adjustment of the epochs of each patch point according to the exposition detailed in Section~\ref{subsec:isochronous_correspondence} and the transformation of the spacecraft states described in Section~\ref{subsec:coord_transf}, once a simulation start epoch has been selected. Although the epochs of the patch points, $\tau_i$, and time intervals, $T_i$, are not part of the design vector under the fixed-time formulation pursued in this work, they are evidently necessary for the propagation of the dynamics along each trajectory segment.

\section{Numerical results}
\label{sec:results}
To evaluate how the proposed LM implementation fares at transitioning CR3BP trajectories into the HFEM, a series of application cases of varying complexity was pondered and numerically simulated. The application of the proposed algorithm is compared with the baseline MN update equation in terms of convergence capabilities and proximity of the converged solution to the initial CR3BP guess. To this end, let $\norm{\mathbf{E}_k}_\mathbf{P}$ denote the position error from the design vector at iteration $k$ to the CR3BP trajectory, i.e. $\mathbf{P}$ is a matrix whose diagonal repeats $(1,1,1,0,0,0)$ for all state entries of the design vector $\mathbf{X}$, with zeros elsewhere.

The numerical simulations are carried out in a Python implementation, with the dynamics being propagated through the TU Delft Astrodynamics Toolbox (TUdat) \cite{tudatspace}, a comprehensive astrodynamics library that computes planetary ephemeris via SPICE kernels. The dynamics are integrated using an adaptive Runge--Kutta--Verner method of order 8(9) with tolerance set to $10^{-12}$. Internal integration of the state transition matrix, for the computation of the Jacobian in Eq.~\eqref{eq:jac}, is performed through a standard adaptive Runge--Kutta method of order 8 with a tolerance set to the same value. Periodic orbits in the CR3BP are retrieved from an extensive NASA catalog \cite{jpl_orbits}, which also provides the relevant characteristic CR3BP scale factors considered in this work. The start epoch for the simulations was set to the 1st of January, 2020. The results obtained will naturally differ for other date choices, but preliminary testing has shown that this effect is not significant for most applications.

In this work, a solution is assumed to have converged once the (dimensionless) continuity residual satisfies $\norm{\mathbf{F}}<10^{-10}$, which corresponds to a deviation of around $39~\si{\milli \meter}/0.1~\si{\micro\meter\per\second}$, if entirely in position/velocity. This degree of accuracy is sufficient for the purposes of this paper, i.e. to analyze and compare the adequacy of the transition algorithms under analysis. Moreover, residuals below this threshold would begin to approach the accuracy of the numerical integration of the dynamics.

An open-source implementation of the algorithm proposed, including the datasets for the test cases to be presented, is available at \url{https://github.com/antoniownunes/LM_mwe}. The reader is encouraged to consult this reference for a better understanding of the inner workings of the algorithms employed and the interface with the TUdat package.

\subsection{Quasi-periodic L2 Northern Halo orbit}
\label{subsec:qpo_L2_halo}
\begin{figure}[h]
	\centering
	\begin{subfigure}{.48\linewidth}
		\centering
		\includegraphics[width=\textwidth,  trim={0 0.5cm 0 2cm}, clip]{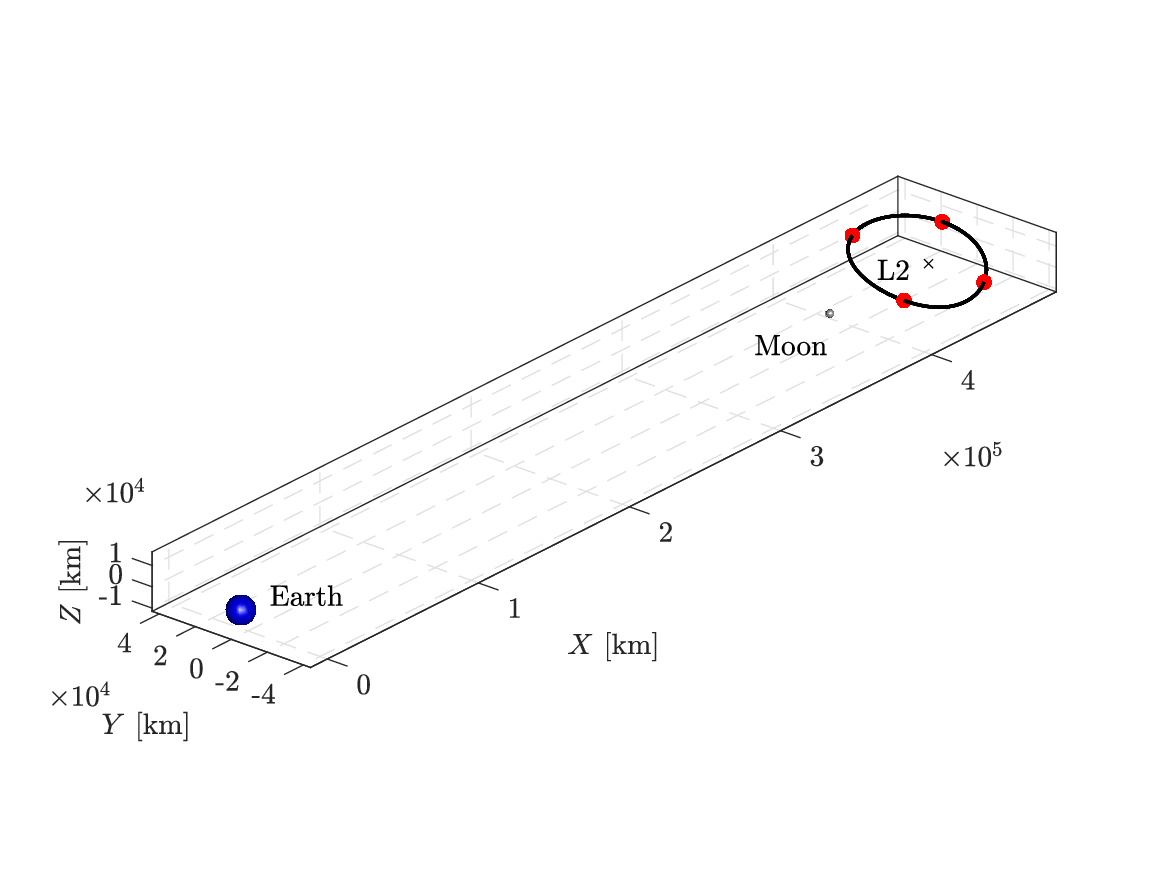}
		\caption{L2 Northern Halo orbit, relative to the Earth and Moon.}
	\end{subfigure}
	\hfill
	\begin{subfigure}{.48\textwidth}
		\centering
		\includegraphics[width=\textwidth,  trim={0 0.5cm 0 2cm}, clip]{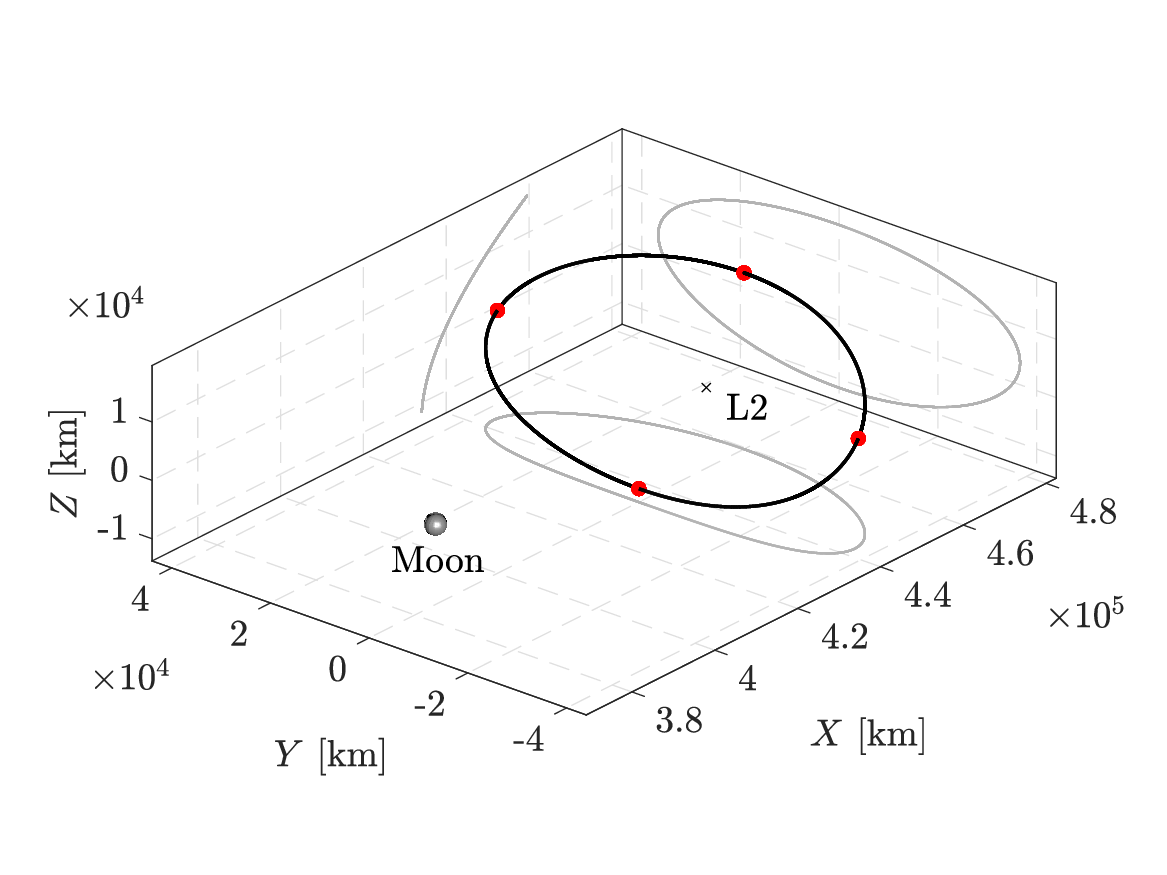}
		\caption{Closeup view of the L2 Northern Halo orbit, relative to the Moon.}
	\end{subfigure}
	\caption{L2 Northern Halo orbit from the CR3BP considered for transition into the HFEM. Dimensions to scale, including Earth and Moon sizes.}
	\label{fig:L2_Halo}
\end{figure}

As a first test case, the transition of a periodic orbit in the CR3BP was considered. This orbit, represented in Fig.~\ref{fig:L2_Halo}, belongs to the Northern Halo family of the L2 Lagrange equilibrium point of the CR3BP dynamics, with a period of approximately $15.08~\text{days}$ and $Z$-amplitude of $2.79\times10^4~\si{\kilo\meter}$. Following the study presented in \cite{park2025CharacterizationL2Analogs}, this particular orbit is characterized by a period that admits the existence of quasi-periodic HFEM solutions in its proximity. Four patch points equally spaced in time were initially selected, as highlighted in the provided figure, corresponding to apolune, perilune, and two symmetrically-placed intermediate points. The objective of this application is to obtain a quasi-periodic solution in proximity to the CR3BP trajectory over a window of time spanning 50 revolutions, or around $2$ years. To do so, the transition algorithms are fed with an initial guess that stacks the four aforementioned patch points 50 times. Although modest, the target duration of 2 years is sufficient to illustrate key differences between the proposed LM implementation and the MN baseline, while also remaining computationally tractable.

\begin{figure}[h]
	\centering
	\begin{minipage}[c]{.48\linewidth}
		\begin{subfigure}[c]{\linewidth}
			\centering
			\includegraphics[width=.87\textwidth]{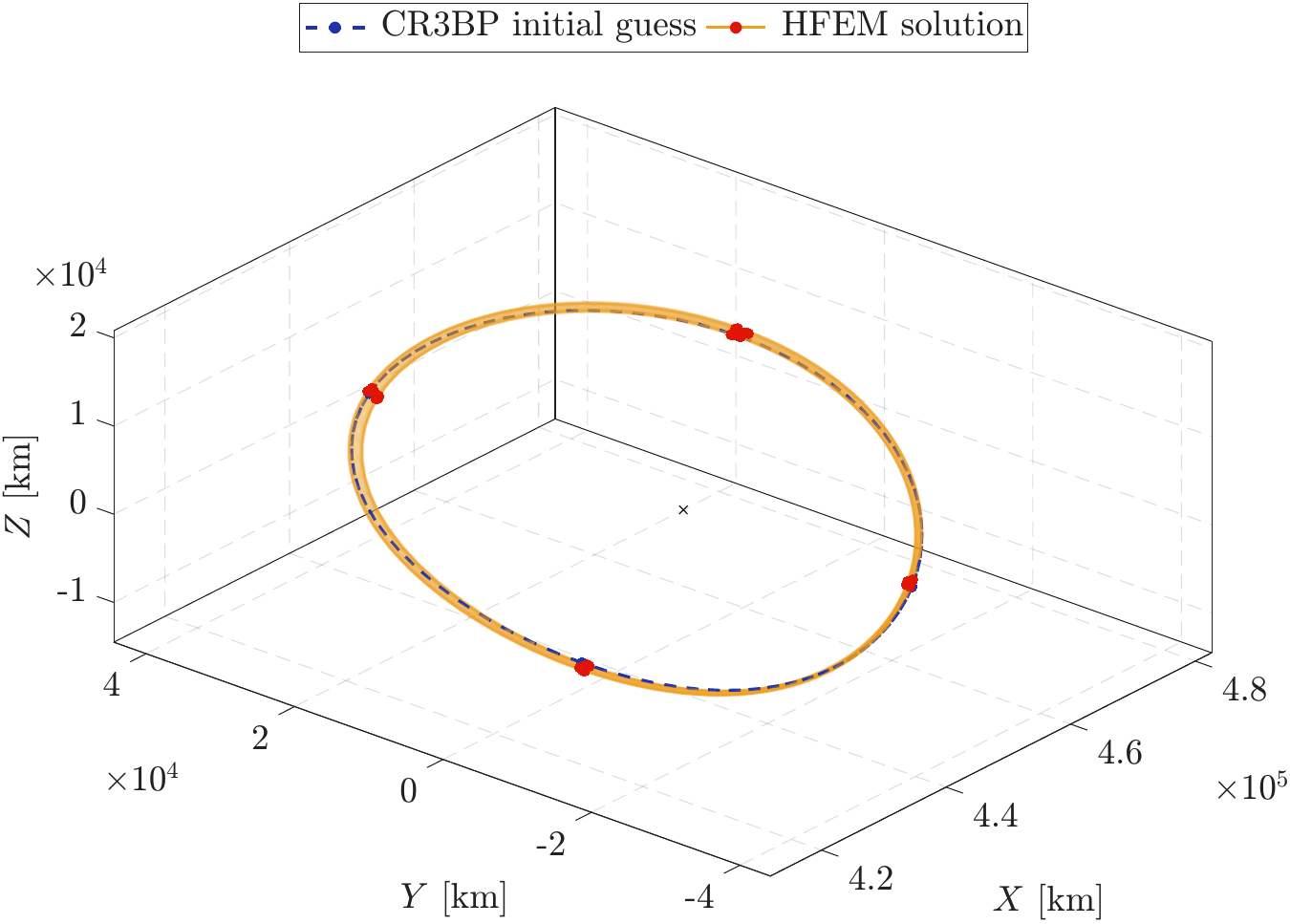}
			\caption{Converged solution via LM algorithm with no proximity constraints.}
			\label{fig:L2_Halo_4pp_sol_LM}
		\end{subfigure}
		\begin{subfigure}[c]{\linewidth}
			\centering
			\includegraphics[width=.87\textwidth]{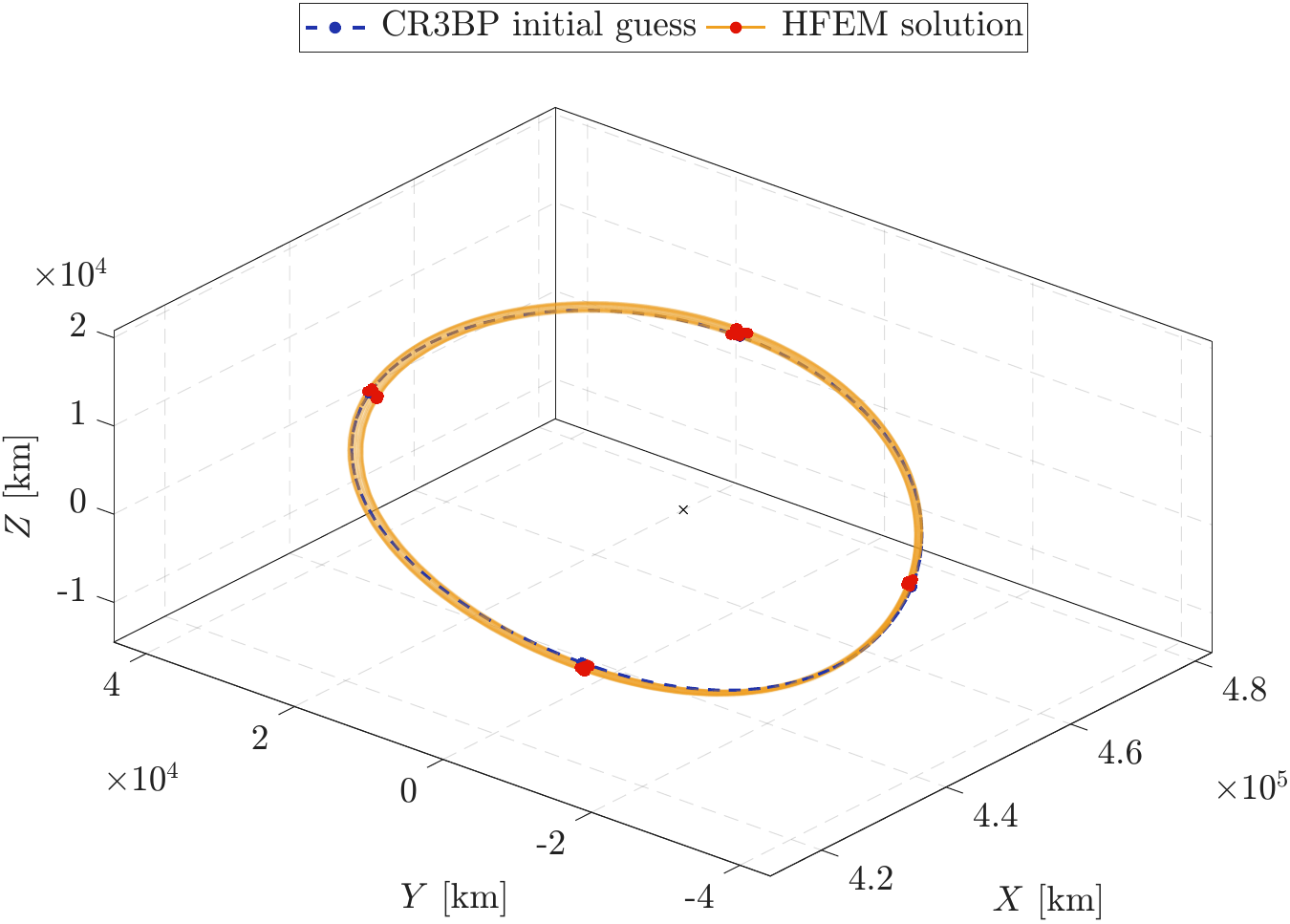}
			\caption{Converged solution via baseline MN update equation.}
			\label{fig:L2_Halo_4pp_sol_MN}
		\end{subfigure}
	\end{minipage}
	\hfill
	\begin{minipage}[c]{.48\linewidth}
		\begin{subfigure}[c]{\linewidth}
			\centering
			\includegraphics[width=\textwidth]{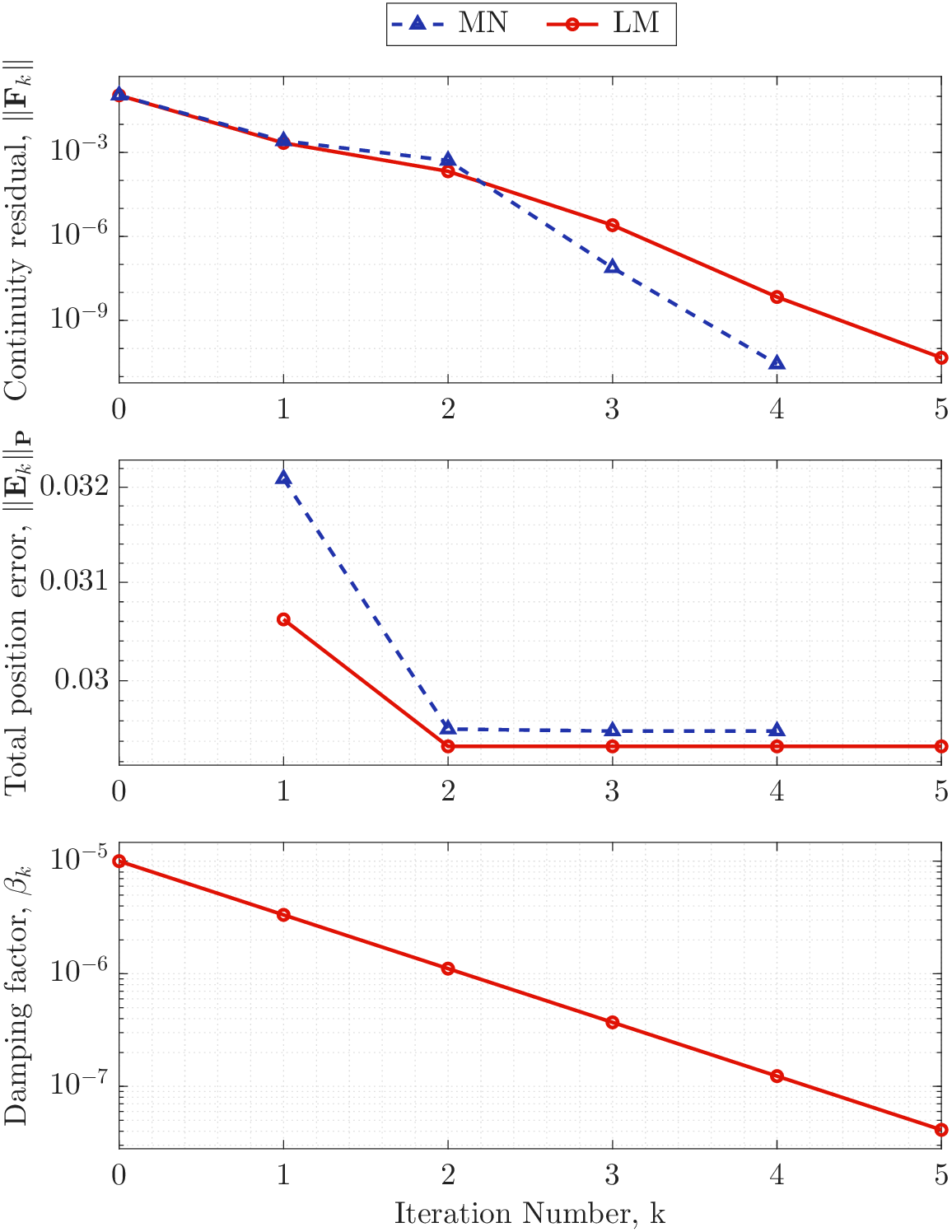}
			\caption{Evolution of meaningful parameters stemming from the application of the transition algorithms.}
			\label{fig:L2_Halo_4pp_res}
		\end{subfigure}
	\end{minipage}
	\caption{Application of the algorithms to transition an L2 Northern Halo orbit from the CR3BP into a 50 revolution quasi-periodic counterpart in the HFEM, considering four patch points per revolution.}
\end{figure}

Given that we expect to find a quasi-periodic solution in the vicinity of the original CR3BP orbit, one may ponder initially how the LM algorithm performs considering no proximity constraints, i.e. setting $\mathbf{Q}=\mathbf{0}$. This also allows for a more direct comparison with the baseline MN implementation. The converged solution resulting from the application of the proposed LM algorithm is provided in Fig.~\ref{fig:L2_Halo_4pp_sol_LM} against the initial CR3BP guess, in the synodic reference frame. Similarly, the converged solution stemming from the  application of the MN update equation is provided in Fig.~\ref{fig:L2_Halo_4pp_sol_MN}. In Fig.~\ref{fig:L2_Halo_4pp_res}, the dimensionless evolution of the continuity residual, $\norm{\mathbf{F}_k}$, and position error, $\norm{\mathbf{E}_k}_\mathbf{P}$, from the application of both algorithms, is plotted in a logarithmic scale. Note that the initial position error is zero and hence cannot be plotted in the logarithmic representation. In addition, the evolution of the damping factor for the LM algorithm is also provided in Fig.~\ref{fig:L2_Halo_4pp_res}. The initial value of $\beta_0=10^{-5}$ was selected after preliminary trials on the quality of the linear approximation, though the effect of changing this parameter will shortly be assessed. Moreover, the adjustment factors used internally by the LM algorithm to adjust $\beta_k$ between iterations,  according to Algorithm~\ref{alg:LM}, were selected as $(\alpha,\eta)=(0.33,2)$. Further discussion on the the pertinence of this selection is also avoided for now, being left for Section~\ref{subsec:more_orbits}, once application of the algorithm is generalized to other trajectories.

By the analysis of Fig.~\ref{fig:L2_Halo_4pp_sol_LM}, one sees that the quasi-periodic solution retrieved by the proposed algorithm is in good agreement with the initial guess, as theoretically postulated. Note that this is true even though the proximity constraint has been disregarded, which attests to a good correspondence between the CR3BP and HFEM models in this particular test case. In fact, there is also a striking similarity, at least qualitatively speaking, between the converged solution stemming from the LM algorithm and the MN update equation, in Fig.~\ref{fig:L2_Halo_4pp_sol_MN}. In quantitative terms, the two approaches generate a solution that resembles the initial patch points to a similar degree which, assuming $\norm{\mathbf{E}}_\mathbf{P}$ has a uniform distribution, corresponds to a position error at each patch point of around $800~\si{\kilo \meter}$ or $3\%$ of the orbit's $Z$-amplitude. The evolution of $\beta_k$ highlights the quality of the linear approximation during the application of the LM algorithm, as a monotonous decrease over the iterations is seen. Additional testing has shown that changing $\beta_0$ had nearly no effect on the final converged solution for this particular example, being only responsible for changes in the convergence speed. For the case presented, the results in Fig.~\ref{fig:L2_Halo_4pp_res} show that the proposed LM scheme and the baseline MN alternative converge almost equally as fast towards a solution that meets the continuity constraint up to a desirable threshold. This is directly reflected in the algorithm run time, considering how the MN and LM implementations take virtually the same time to compute each iteration, assuming the linear approximation is deemed adequate, which is the case.

\begin{figure}[h]
	\centering
	\begin{minipage}[c]{.48\linewidth}
		\begin{subfigure}[c]{\linewidth}
			\centering
			\includegraphics[width=.87\textwidth]{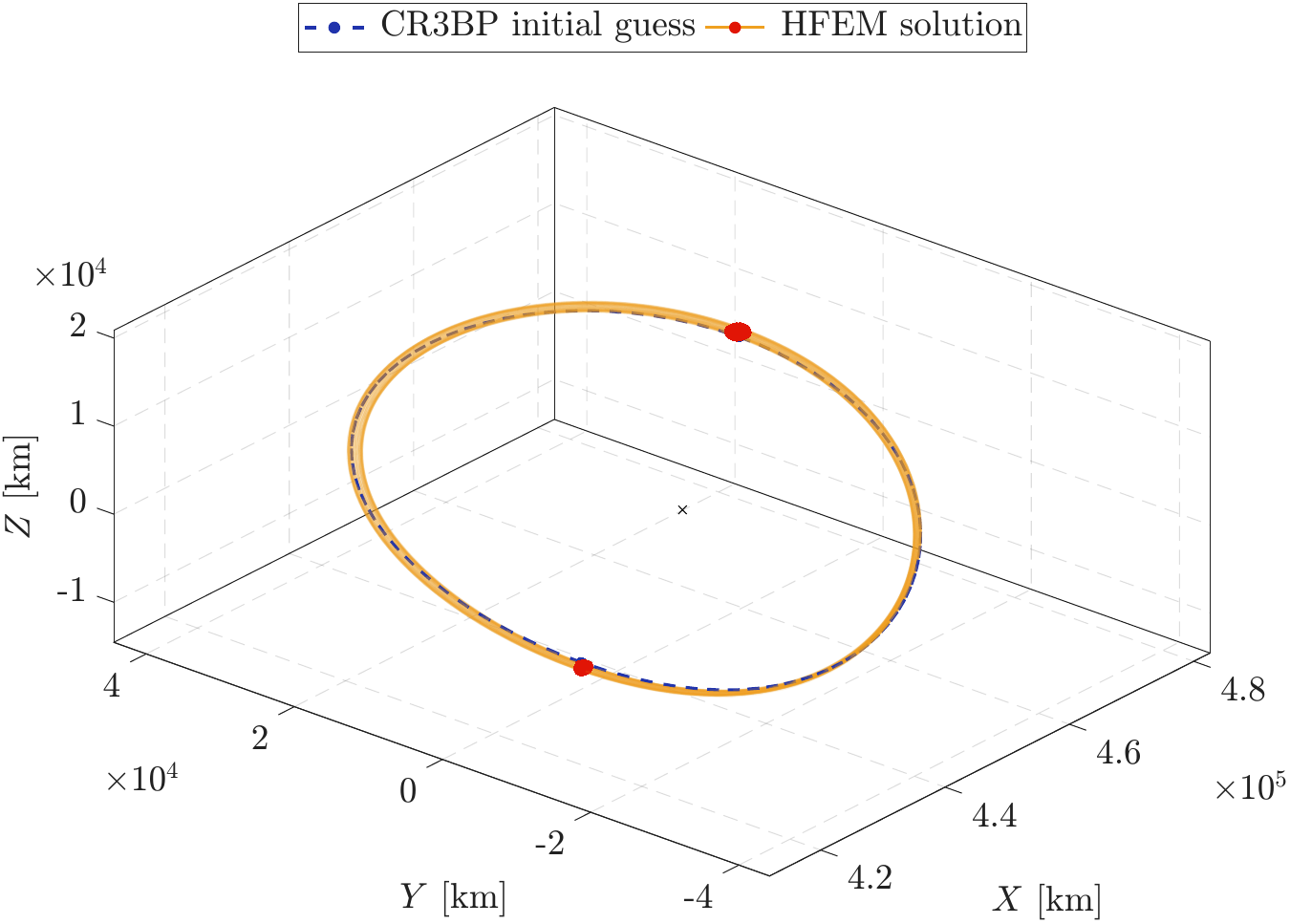}
			\caption{Converged solution via LM algorithm with no proximity constraints, for $\beta_0=0.1$.}
			\label{fig:L2_Halo_2pp_sol_LM_1e-1}
		\end{subfigure}
		\begin{subfigure}[c]{\linewidth}
			\centering
			\includegraphics[width=.87\textwidth]{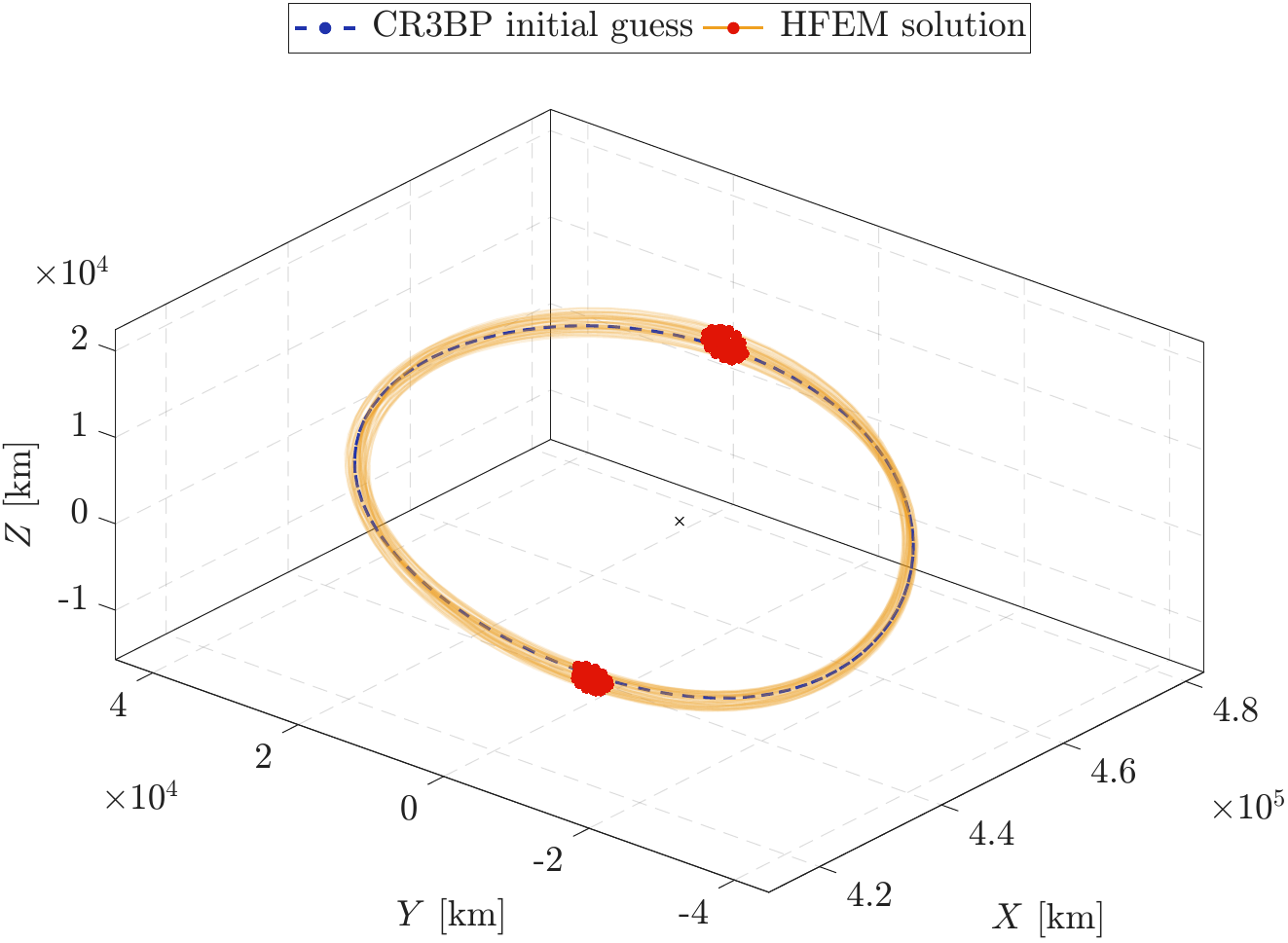}
			\caption{Converged solution via baseline MN update equation.}
			\label{fig:L2_Halo_2pp_sol_MN}
		\end{subfigure}
	\end{minipage}
	\hfill
	\begin{minipage}[c]{.48\linewidth}
		\begin{subfigure}[c]{\linewidth}
			\centering
			\includegraphics[width=\textwidth]{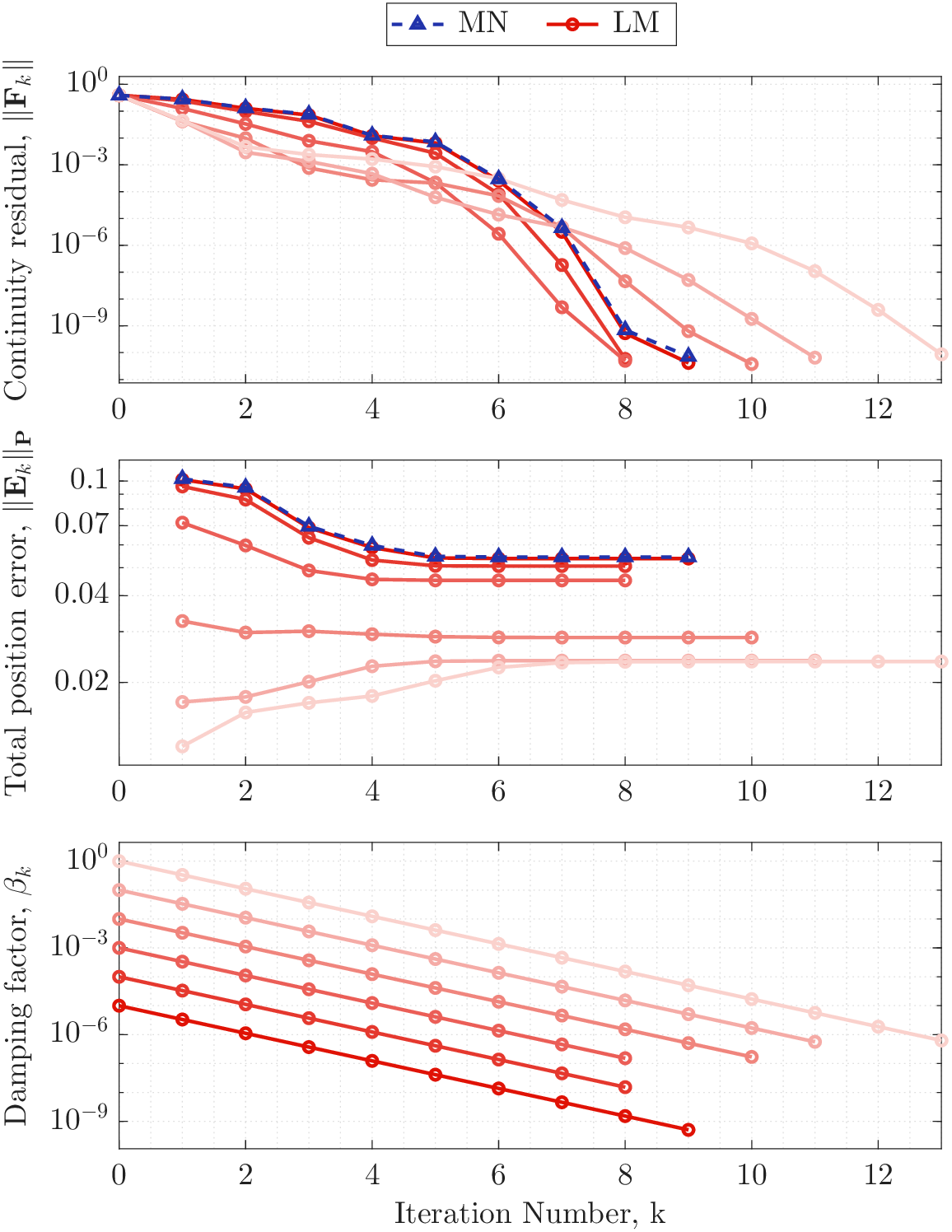}
			\caption{Evolution of meaningful parameters stemming from the application of the transition algorithms. LM algorithm evaluated under different initial values of the damping factor, evidenced in different shades according to the $\beta_k$ plot.}
			\label{fig:L2_Halo_2pp_res}
		\end{subfigure}
	\end{minipage}
	\caption{Application of the algorithms to transition an L2 Northern Halo orbit from the CR3BP into a 50 revolution quasi-periodic counterpart in the HFEM, considering two patch points per revolution.}
	\label{fig:L2_Halo_2pp}
\end{figure}

To further test the capabilities of the LM algorithm, a case where only 2 patch points per revolution are considered was also evaluated. Recalling Fig.~\ref{fig:L2_Halo}, the points at perilune and apolune were selected for this purpose, i.e. the points closest and farthest from the Moon, though additional analysis here omitted has shown that other pair selections yield similar results. In this case, the initial CR3BP guess is found to be less adequate and thus the damping factor for the LM approach plays a particularly important role, given how it leverages the level of trust put on this approximation. To assess this, we evaluate the effect of changing the initial value of the damping factor. As such, Fig.~\ref{fig:L2_Halo_2pp_res} presents the evolution of the continuity residual, position error, and damping factor for various choices of $\beta_0$, evidenced through different shades according to the information in the $\beta_k$ plot. These curves are compared to the baseline behavior of the MN scheme.
Fig.~\ref{fig:L2_Halo_2pp_sol_LM_1e-1} presents the converged solution stemming from the application of the LM algorithm without proximity constraints, considering $\beta_0=0.1$, which is put against the trajectory resulting from the MN baseline, in Fig.~\ref{fig:L2_Halo_2pp_sol_MN}.

The analysis of Fig.~\ref{fig:L2_Halo_2pp} evidences that reducing the number of patch points per revolution might impact the quality of the trajectory retrieved through the MN approach, though this effect may be attenuated under the LM algorithm through an appropriate selection of $\beta_0$. In particular, we find that the initial trust put on the linear approximation has a relevant influence on which solution the LM algorithm converges to. This leads to the important conclusion that, while the linear approximation may be deemed adequate from the standpoint of reducing the continuity residual, it may still not provide a reasonable approximation of the intricate dynamics at play. In fact, we see that as $\beta_0$ is decreased, i.e. we trust the linear approximation and allow the initial update steps to be larger, the LM approach tends to the behavior of the MN update equation, reaching convergence to a solution of similar proximity to the CR3BP guess in exactly the same number of iterations. We note, however, that a more appropriate choice of initial damping factor leads to converged solutions that retain better proximity with respect to the CR3BP trajectory, which may be observed qualitatively through the comparison of Fig.~\ref{fig:L2_Halo_2pp_sol_LM_1e-1} and Fig.~\ref{fig:L2_Halo_2pp_sol_MN}, where the solutions stemming from the LM algorithm with $\beta_0=0.1$ and MN update equation, respectively, are provided. This is further confirmed in quantitative terms by the position error, which is halved with respect to the MN solution for such a choice of $\beta_0$, as shown in Fig.~\ref{fig:L2_Halo_2pp_res}. Evidently, as $\beta_0$ is increased, convergence is delayed due to the strong penalization on the norm of the update step, and reduced benefits are observed when increasing $\beta_0$ past $10^{-1}$. Nonetheless, these considerations show that, even as a base strategy with no proximity constraints, the proposed LM algorithm has more flexibility than the MN counterpart when it comes to addressing the concerns related to the linear approximation at the basis of convexification, simply by enacting changes to the initial damping factor. To this end, one may freely tune the proposed algorithm to match the MN update equation or to take slightly longer to achieve a more desirable alternative. In fact, given that the converged solution considering only two patch points may be brought so close to the original CR3BP orbit over 50 revolutions with a proper choice of $\beta_0$, achieving similar results to the previous case with 4 patch points, the inclusion of the proximity constraint was deemed not relevant for study under this example.

\subsection{Additional quasi-periodic trajectories}
\label{subsec:more_orbits}
\begin{figure}[h]
	\centering
	\includegraphics[width=\textwidth]{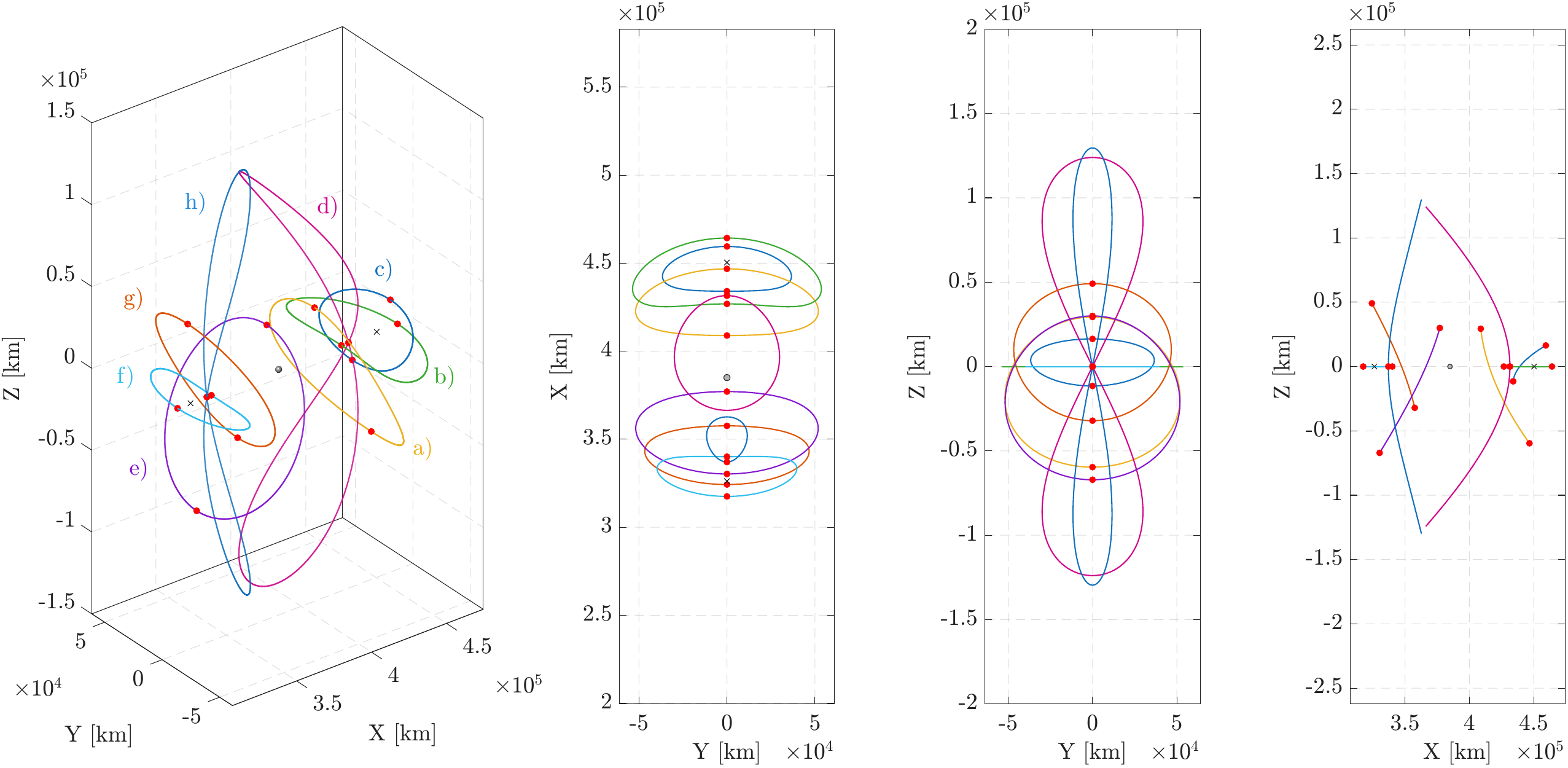}
	\caption{Additional CR3BP periodic orbits considered for transition into the HFEM. The trajectories belong to the Southern Halo, Lyapunov, Northern Halo, and Vertical families about the L2 and L1 Lagrange points and are identified by a) through h), in that order. Dimensions to scale, including Moon size.}
	\label{fig:more_orbits}
\end{figure}

In order to generalize the findings of the previous section, we evaluate the performance of the proposed LM algorithm and MN baseline over a wide range of orbits with different periods and belonging to various families. Focusing on the L1 and L2 Lagrange points, we study the transition of CR3BP periodic orbits from the Southern Halo, Lyapunov, Northern Halo, and Vertical families. These orbits are displayed in Fig.~\ref{fig:more_orbits} and are discerned by their color and an identifier, namely a) through h). Besides the trajectories in three-dimensional space, their projections along each plane of the synodic reference frame are also provided for a more complete visualization. The information regarding orbital family and period is contained in Table~\ref{tab:more_orbits_periods}. Note that the trajectory studied in Section~\ref{subsec:qpo_L2_halo} is here recovered and corresponds to orbit c). Similarly to the previous test case, these orbits were carefully chosen so as to avoid transition-challenging regions, which are out of the scope of this paper.

\begin{table}[h]
	\centering
	\caption{Periods and families of the CR3BP orbits considered in Fig.~\ref{fig:more_orbits}.}
	\label{tab:more_orbits_periods}
	\begin{tabular}{llll}
		\hline
		Lagrange Point  & Orbit & Family & Period [days]\\
		\hline
		\multirow{4}{*}{L2}& a) & Southern Halo & 14.02\\
		& b) & Lyapunov & 17.02\\
		& c) & Northern Halo & 15.08\\
		& d) & Vertical & 22.09\\
		\hline
		\multirow{4}{*}{L1}& e) & Southern Halo & 11.50\\
		& f) & Lyapunov & 12.72\\
		& g) & Northern Halo & 12.33 \\
		& h) & Vertical & 21.72 \\
		\hline
	\end{tabular}
\end{table}

For the sake of tractability, transition is once more targeted over 50 revolutions for each trajectory. As highlighted in Fig.~\ref{fig:more_orbits}, the shooting problem is posed considering two patch points per revolution of the original CR3BP structures, corresponding to apolune and perilune --- in the case of the Vertical orbits, d) and h), both points coincide in position. Aiming for completeness, the MN baseline is evaluated considering different thresholds for the update norm, $\gamma$, according to Eq.~\eqref{eq:MN_attenuation}: $\gamma=\infty$ for an unbounded application of the method, $\gamma=0.1$ to reflect typical choices in literature, and $\gamma=0.01$ to potentially remain closer to the original CR3BP orbit. This is compared against three choices of $\beta_0$ representing decreasing levels of initial trust on the linear approximation, namely $\beta_0=10^{-5}$, $\beta_0=10^{-3}$, and $\beta_0=10^{-1}$. Much like in Section~\ref{subsec:qpo_L2_halo}, the adjustment factors used internally by the LM algorithm were selected as $(\alpha,\eta)=(0.33,2)$. In \ref{ap:param_study}, algorithmic performance following common strategies for the selection of these factors is presented in the form of a parametric study, applied over the orbits in Table~\ref{tab:more_orbits_periods}. The results therein attest to the aforementioned selection when prioritizing proximity to the original CR3BP structures, which is evidently relevant in this application case. As also approached in \ref{ap:param_study}, other choices may yield benefits in terms of speed or adaptability, boosting the algorithm's capacity of handling challenging transitions at the cost of the proximity measure. This is later exploited in Section~\ref{subsec:transfer}.

The results stemming from the application of the baseline MN update equation and the proposed LM algorithm are presented in Table~\ref{tab:more_orbits_results}, namely the number of iterations necessary to reach convergence, denoted $m$, and the final total position error of the converged solution, $\norm{\mathbf{E}_m}_\mathbf{P}$. In the cases where MN diverges, i.e. the algorithm departs significantly from the initial conditions, we indicate this with ``Div.''. When the algorithms are unable to converge within $100$ iterations to the established threshold of $\norm{\mathbf{F}(\mathbf{X})}<10^{-10}$, their execution is stopped and the data in the table is absent. This could represent slow convergence or stalling at a local minimum that does not respect the continuity constraint to the degree required.

\begin{table}[h]
	\centering
	\caption{Results from the application of the algorithms to transition the CR3BP orbits of Fig.~\ref{fig:more_orbits}, for various values of $\gamma$, for the MN baseline, and $\beta_0$, for the proposed LM algorithm.}
	\label{tab:more_orbits_results}
	\makebox[\linewidth][c]{%
	\begin{tabular}{lllllllllllll}
		\hline
		\multirow[b]{3}{*}{Orbit} & \multicolumn{6}{l}{Number of iterations, $m$} & \multicolumn{6}{l}{Final total position error, $\norm{\mathbf{E}_m}_\mathbf{P}~\left( \times 10^{-2}\right)$}\\
		\cmidrule(lr){2-7} \cmidrule(lr){8-13}
		& \multicolumn{3}{l}{MN ($\gamma$)} & \multicolumn{3}{l}{LM ($\beta_0$)} & \multicolumn{3}{l}{MN ($\gamma$)}  & \multicolumn{3}{l}{LM ($\beta_0$)}\\
		\cmidrule(lr){2-4} \cmidrule(lr){5-7} \cmidrule(lr){8-10} \cmidrule(lr){11-13}
		& $\infty$ & $10^{-1}$ & $10^{-2}$ & $10^{-5}$ & $10^{-3}$ & $10^{-1}$ & $\infty$ & $10^{-1}$ & $10^{-2}$  & $10^{-5}$ & $10^{-3}$ & $10^{-1}$ \\
		\hline
		a) & Div. & 26 & -- &  15 &  16 &  56 & Div. & 24.95 & -- &  6.520 &  7.817 &  3.095 \\
		b) & Div. & -- & -- &  14 &  12 &  16 & Div. & -- & -- &  2.994 &  5.841 &  2.338 \\
		c) & 9 & 7 & 17 &  9 &  8 &  11 & 5.447 & 5.096 & 3.723 &  5.383 &  4.519 &  2.379 \\
		d) & Div. & Div. & -- &  -- &  14 &  18 & Div. & Div. & -- &  -- &  5.127 &  4.832 \\
		\hline
		e) & Div. & 17 & -- &  12 &  23 &  -- & Div. & 12.04 & -- &  7.973 &  7.416 &  -- \\
		f) & Div. & 10 & -- &  9 &  10 &  13 & Div. & 15.99 & -- &  5.488 &  4.417 &  1.224 \\
		g) & 9 & 9 & 46 &  8 &  9 &  13 & 6.744 & 6.044 & 5.670 &  5.948 &  3.949 &  2.643\\
		h) & Div. & Div. & -- &  26 &  6 &  17 & Div. & Div. & -- &  16.05 &  14.88 &  4.948 \\
		\hline
	\end{tabular}
	}
\end{table}

\begin{figure}[h]
	\centering
	\includegraphics[width=0.98\textwidth]{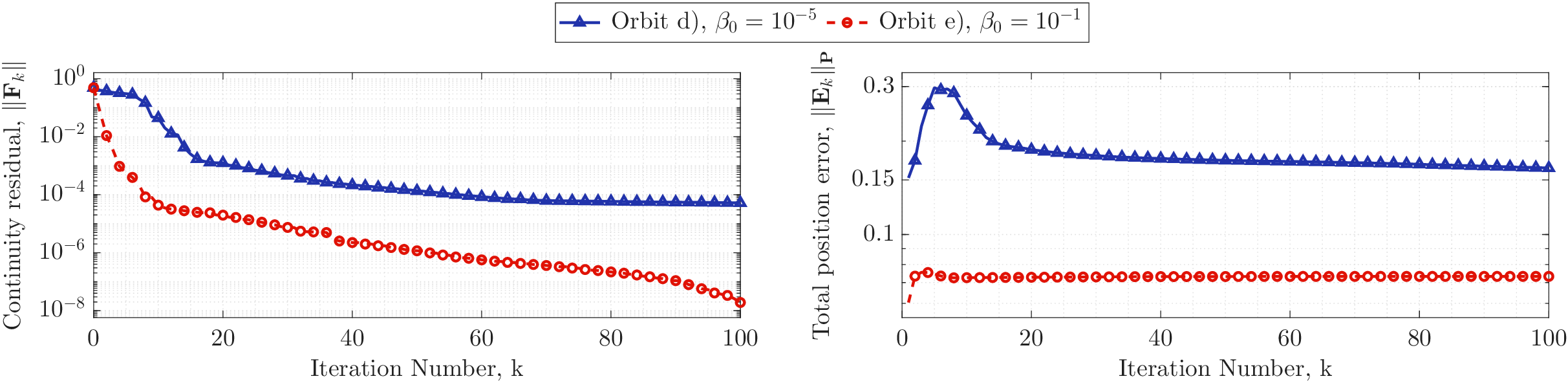}
	\caption{Evolution of the continuity residual and total position error for the two test cases in Table~\ref{tab:more_orbits_results} for which the LM algorithm is unsuccessful.}
	\label{fig:additional_results}
\end{figure}

The results in Table~\ref{tab:more_orbits_results} underscore key differences between the two algorithms and generalize the findings of Section~\ref{subsec:qpo_L2_halo}. To this end, perhaps the most relevant observation is that the MN baseline is susceptible to divergence, especially if the norm of the update step is not bounded, as evidenced by the fact that convergence is only achieved for two of the eight CR3BP orbits considered when $\gamma=\infty$. For the most part, this may be corrected by limiting the norm of the update step through $\gamma$, though two concerns arise. On the one hand, the optimal choice for $\gamma$ is not necessarily evident, since there is no prescribed value that works in all cases, as the table goes to show. On the other hand, it may also be the case that the attenuation necessary to avoid divergence significantly slows down the algorithm, which results in cases where convergence was not possible for any $\gamma$, at least within 100 iterations --- such as with orbits b), d), and h). Evidently, it is not possible to test all values of $\gamma$, and it may be the case that a very particular choice is successful over these challenging cases. Similarly, such issues may potentially be resolved by increasing the number of patch points per revolution, changing their location, or switching to a more complete formulation of the shooting problem (such as variable-time). Nonetheless, the results presented evidence that the MN update scheme may suffer in terms of robustness in specific applications, especially when given poor initial guesses --- here emulated by limiting the number of patch points to two per revolution. In stark contrast, the inherent stability of the proposed LM algorithm stands out for the test cases considered, as divergence concerns are completely mitigated. In fact, only for orbit d) with $\beta_0=10^{-5}$ and orbit e) with $\beta_0=10^{-1}$ does the LM algorithm struggle to converge within 100 iterations. To further analyze these two cases, Fig.~\ref{fig:additional_results} presents the corresponding evolution of the continuity residual and position error, with markers at every second iteration.

Fig.~\ref{fig:additional_results} highlights that the apparent failure for orbit e) with $\beta_0=10^{-1}$~is, in reality, a case of slow convergence, prompted by the relatively large initial damping factor --- at the 100th iteration, $\norm{\mathbf{F}}\approx1.92\times10^{-8}$. On the contrary, convergence for orbit d) with $\beta_0=10^{-5}$ appears to falter due to the LM algorithm initially straying away from the CR3BP guess, likely because of the high level of trust put on the linear approximation at the start. In \ref{ap:param_study}, we show that difficulties in transitioning this orbit with small initial damping are transversal to other choices of the adjustment factors $(\alpha,\eta)$, confirming a significant lack of quality of the linear approximation in this specific case. This sheds light on the failure for $\beta_0=10^{-5}$ and~motivates the choice of larger values, according to the previous results in Table~\ref{tab:more_orbits_results}.~These examples further illustrate that, while the initial damping factor may be selected freely without a consequence to the numerical stability of the scheme, a poor selection may negatively influence exactly \textit{how} the algorithm evolves.

Returning the attention to Table~\ref{tab:more_orbits_results}, we also remark that, at least for the test cases considered, the proposed LM algorithm may always be made to converge to a final solution that is closer to the original CR3BP orbit through an appropriate selection of $\beta_0$, as evidenced by the final total position error. This is true even though the MN baseline actively minimizes the norm of each update step between iterations, precisely to remain close to the initial guess, which the LM algorithm does not. And, while decreasing $\gamma$ may successfully decrease $\norm{\mathbf{E}_m}_\mathbf{P}$, we still see that the values achieved by the MN baseline are at least $50\%$ larger than the best LM results. To this end, in line with the discussion in Section~\ref{subsec:qpo_L2_halo}, it is typically the case that increasing $\beta_0$ leads to a final solution that is closer to the initial CR3BP guess, at the cost of requiring more iterations to reach convergence. However, some exceptions may be found. Considering also the failure cases discussed above, these exceptions suggest that it may be necessary to tune $\beta_0$ empirically and on an application basis, much like the MN scheme with $\gamma$, which may be seen as a drawback shared by both approaches. This information should be leveraged with the improvements in robustness of the proposed scheme, in comparison with the MN alternative, when choosing the ideal strategy to tackle a specific transition problem.

\subsection{Orbital transfer}
\label{subsec:transfer}
To further evaluate the proposed LM algorithm, a test case concerning the transition of an orbital transfer between two quasi-periodic trajectories was considered. This is a topic that has gained particular traction in recent research preceding and following the successful maneuver performed by NASA's Artemis P1 spacecraft between L2 and L1 quasi-planar Lissajous trajectories \cite{broschart2009Artemis}. In pursuit of a similar test scenario, a transfer from the L2 Northern Halo orbit to the L1 Lyapunov planar orbit, from Table~\ref{tab:more_orbits_periods}, was pondered. This application differs slightly from the transfer pursued by the Artemis P1 spacecraft in that we consider no efforts to firstly bring the spacecraft to a quasi-planar trajectory before the transfer occurs. In fact, we note that, due to the significant spatial differences between the L2 Northern Halo and L1 Lyapunov orbits, a solution that closely matches the two CR3BP trajectories is unlikely to be reached without the introduction of a control maneuver. Nonetheless this application scenario is deemed interesting from the perspective of putting to the test the capabilities of the proposed algorithm and to verify how the final solution adapts to overcome the additional challenge of departing a three-dimensional orbit to reach a planar counterpart. To accommodate for the spatial discrepancy between the two orbits, several revolutions before and after the transfer will be considered. 

To ensure that the transfer occurs without a change to the spacecraft's energy, the two orbits were carefully chosen for their correspondence in terms of the Jacobi integral\footnote{Although this quantity is not conserved under the non-autonomous HFEM, which admits no constants of motion, the choice of orbits that respect its conservation in the CR3BP is expected to provide better guarantees for the feasibility of a transfer trajectory between them under the high-fidelity dynamics.}, which takes the value of $\mathcal{J}\approx3.1445$. The initial guess for the transfer trajectory is constructed from a pair of unstable--stable manifolds that depart the L2 Northern Halo orbit and arrive at the L1 Lyapunov orbit, respectively, and which near intersection in an intermediate region near the Moon. The unstable (stable) manifold was computed by introducing an infinitesimal perturbation at various points along the nominal orbits in the local direction of the unstable (stable) linear mode, retrieved through Floquet theory\footnote{Consult the book in \cite{ross3BPbook} for a great introduction on this topic within the context of the three-body problem.}, and then propagating the dynamics forward (backward) in time. By using nearly-intersecting invariant directions to design an initial~guess for the transfer, we assume that the algorithms naturally converge towards the trajectory of ``least resistance'' that exploits the dynamics at play to meet the transfer goal without the need for control maneuvers. The specific manifolds and corresponding periods of integration to retrieve them are chosen manually to reach a suitable pair for the purposes of the transfer. To this end, no further efforts were put forth to optimize the manifolds to provide the best initial guess. Future work could focus on the computation of heteroclinic connections for this purpose \cite{mccarthy2023heteroclinic}, but such an approach was considered outside of the scope of this study. In fact, one may also interpret the discontinuity at the transition point between the two manifolds as a test to the numerical stability and convergence capabilities of the algorithms employed.

\begin{figure}[h]
	\centering
	\begin{subfigure}{.48\linewidth}
		\centering
		\includegraphics[width=\textwidth]{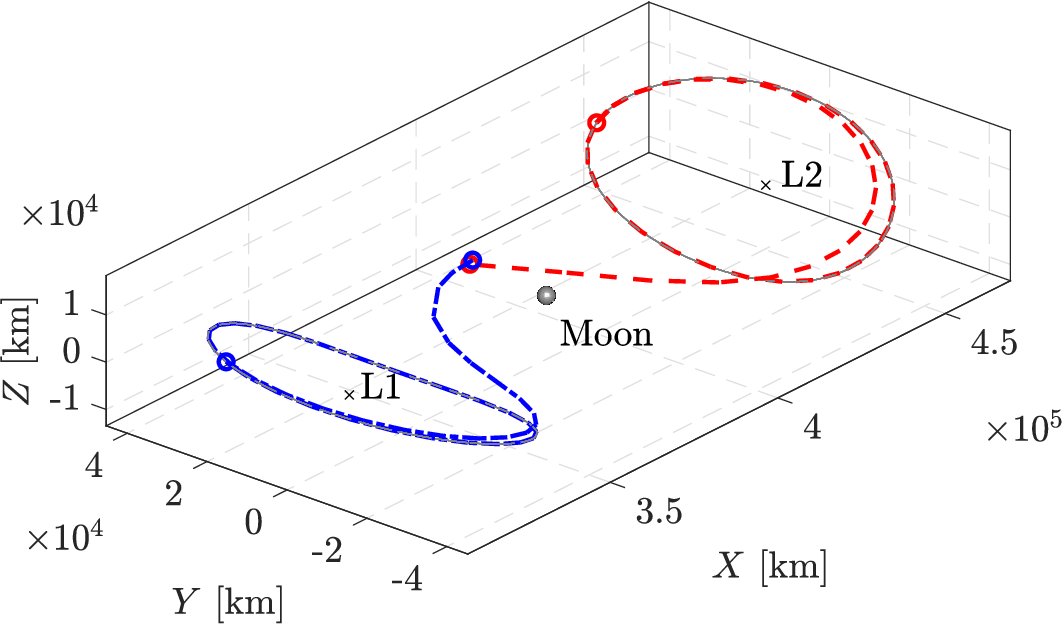}
		\caption{Global view of the orbits and stable-unstable manifold pair for the transfer.}
		\label{fig:transfer_segments}
	\end{subfigure}
	\hfill
	\begin{subfigure}{.48\textwidth}
		\centering
		\includegraphics[width=\textwidth]{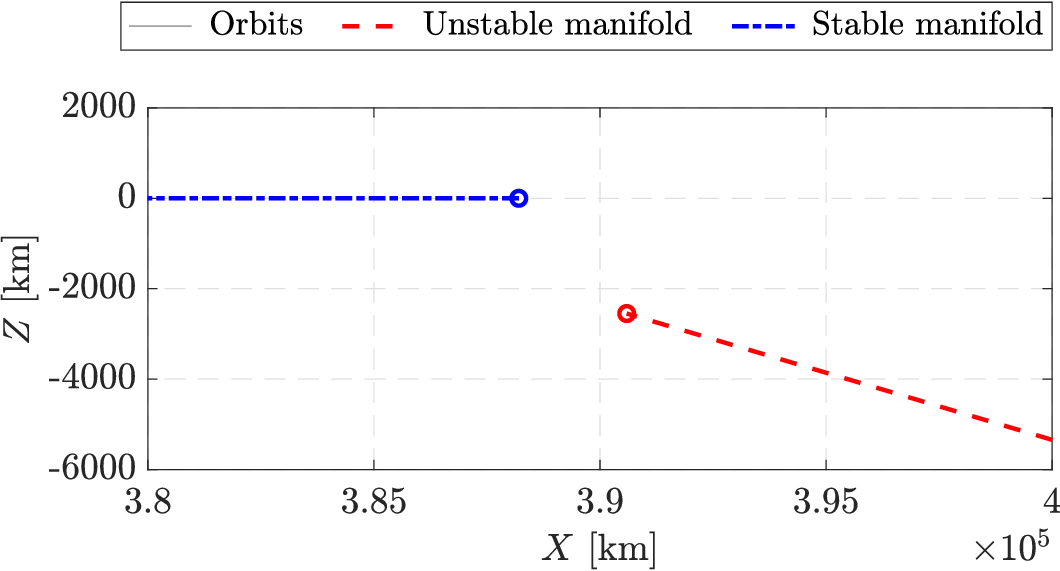}
		\caption{Detailed $X$--$Z$ view of the manifold transition region.}
		\label{fig:transfer_detail}
	\end{subfigure}
	\caption{Construction of the initial CR3BP guess for a transfer trajectory between the L2 Northern Halo orbit and L1 Lyapunov orbit using an unstable manifold from the former and stable manifold from the latter.}
	\label{fig:transfer}
\end{figure}

The CR3BP L2 Northern Halo and L1 Lyapunov orbits are represented in Fig.~\ref{fig:transfer_segments}, alongside the respective unstable and stable manifolds considered, whose endpoints are denoted by a circle. Note that both manifolds revolve once in close proximity to their respective orbit before departing in the second revolution. In Fig.~\ref{fig:transfer_detail}, a detailed $X$--$Z$ view of the region where the transition from the unstable manifold to the stable manifold occurs is provided to evidence the discontinuity in spacecraft position and heading.

The trajectory to be transitioned into the HFEM was built considering the following sections: \textit{(a)} $10$ full revolutions around the L2 Northern Halo orbit, starting at apolune; \textit{(b)} incomplete revolution until the unstable manifold start point; \textit{(c)} departure through the unstable manifold; \textit{(d)} access to the L1 Lyapunov orbit via the stable manifold; \textit{(e)} incomplete revolution until apolune; \textit{(f)} $10$ full revolutions.
In contrast to the isolated orbits studied in Section~\ref{subsec:qpo_L2_halo} and \ref{subsec:more_orbits}, a spacecraft will feel significant changes in velocity as it passes near the Moon during sections \textit{(c)} and \textit{(d)}, and thus partitioning the overall trajectory in equal intervals of time is inadequate. To this end, the solution is initially parameterized by arc-length, i.e. the distance along the trajectory, and then equally spaced segments are chosen instead. To evaluate the robustness of the algorithms under analysis, the number of segments considered (and hence the number of patch points) is varied. For now, no proximity constraint is introduced to better compare their adequacy and convergence capabilities.

The results stemming from the application of the MN update equation and LM algorithm, considering different values of $\gamma$ and $\beta_0$, are shown in Table~\ref{tab:results_transfer} for various choices of the number of segments, $n$. The table highlights the numbers of iterations needed to reach convergence, $m$, according to a continuity residual threshold of $\norm{\mathbf{F}(\mathbf{X})}<10^{-10}$, and the total position error of the final converged solution, $\norm{\mathbf{E}_m}_\mathbf{P}$. As with the analysis presented in Section~\ref{subsec:more_orbits}, the number of iterations is substituted by ``Div.'' in the cases where the MN update equation diverges. Alternatively, if the algorithms are not able to bring the continuity residual below the desired threshold within 100 iterations, their execution is stopped and the corresponding values in the table are absent. In line with the discussion drawn out in \ref{ap:param_study}, the increase in complexity of this application case prompted the selection of $(\alpha,\eta)=(0.2,1.5)$ as the adjustment factors used internally by the LM scheme, following Algorithm~\ref{alg:LM}. This choice contrasts with the one considered in previous sections, since the algorithm's adaptability is here valued over the proximity of the final solution to the original CR3BP structures.


\begin{table}[h]
	\centering
	\caption{Results from the application of the algorithms to transition the CR3BP transfer trajectory between the L2 Northern Halo orbit and L1 Lyapunov orbit, for various numbers of segments, equally spaced in arc-length, considered in the multiple-shooting formulation.}
	\label{tab:results_transfer}
	\makebox[\linewidth][c]{%
		\begin{tabular}{lllllllllllll}
			\hline
			\multirow[b]{3}{2.5cm}{Number of segments, $n$}  & \multicolumn{6}{l}{Number of iterations, $m$} & \multicolumn{6}{l}{Final total position error, $\norm{\mathbf{E}_m}_\mathbf{P}~\left(\times10^{-1}\right)$} \\
			\cmidrule(lr){2-7} \cmidrule(lr){8-13}
			& \multicolumn{3}{l}{MN ($\gamma$)} & \multicolumn{3}{l}{LM ($\beta_0$)}  & \multicolumn{3}{l}{MN ($\gamma$)} & \multicolumn{3}{l}{LM ($\beta_0$)}  \\
			\cmidrule(lr){2-4} \cmidrule(lr){5-7} \cmidrule(lr){8-10} \cmidrule(lr){11-13}
			& $\infty$ & $10^{-1}$ & $10^{-2}$  & $10^{-5}$ & $10^{-3}$ & $10^{-1}$ & $\infty$ & $10^{-1}$ & $10^{-2}$  & $10^{-5}$ & $10^{-3}$ & $10^{-1}$ \\
			\hline
			$70$ & 9 & 18 & 75 &  10 &  16 &   66 & 1.938 & 3.129 & 2.665 &   2.314 &   2.894 &   2.807  \\
			$60$ & 9 & 16 & 86 &   9 &   35 &   81 & 1.929 & 2.726 & 2.604 &   1.800 &   2.377 &   2.539 \\
			$50$ & Div. & 22 & -- &   17 &   17 &   -- & Div. & 2.820 & -- &   2.699 &   2.582 &   -- \\
			$40$ & Div. & -- & -- &   20 &   15 &   94 & Div. & -- & -- &   1.965 &   1.937 &   1.944 \\
			\hline
		\end{tabular}
	}
\end{table}

Table~\ref{tab:results_transfer} serves to further consolidate that, when the initial guess is sufficiently adequate, the LM approach can be made to follow closely the behavior of the baseline MN update equation through an appropriate choice of $\beta_0$, in terms of convergence speed and proximity to the initial CR3BP guess. However, in this application case, we note that the use of a larger initial damping factor does not seem to guarantee better results in terms of the final position error, despite of what the results in the previous sections might have suggested, which might be attributed to the more challenging trajectory profile. In fact, the very slow convergence observed for $\beta_0=10^{-1}$ may be interpreted as evidence for this fact. Nonetheless, the results obtained once more highlight the inherent robustness of the LM algorithm, which is especially apparent when the resolution of the initial CR3BP guess is poor and the unbounded MN alternative diverges, e.g. for $n=50$ and $n=40$. In fact, we see that the LM scheme retains its performance without a sacrifice to the quality of the final converged solution, achieving similar values of $\norm{\mathbf{E}_m}_\mathbf{P}$ for the smaller choices of $\beta_0$. Naturally, as the initial guess becomes increasingly inadequate, both algorithms generally require more iterations to reach a converged solution. Further analysis shows that, for $n<40$, the number of segments is insufficient to achieve a converged solution that brings the continuity residual to the desired threshold within 100 iterations. Even so, the fact that there are specific scenarios where the LM approach is able to converge whereas the MN alternative diverges represents palpable evidence of the benefits in robustness provided by the algorithm proposed in this work.

\begin{figure}[h]
	\centering
	\begin{minipage}[c]{.48\linewidth}
		\begin{subfigure}[c]{\linewidth}
			\centering
			\includegraphics[width=.83\textwidth]{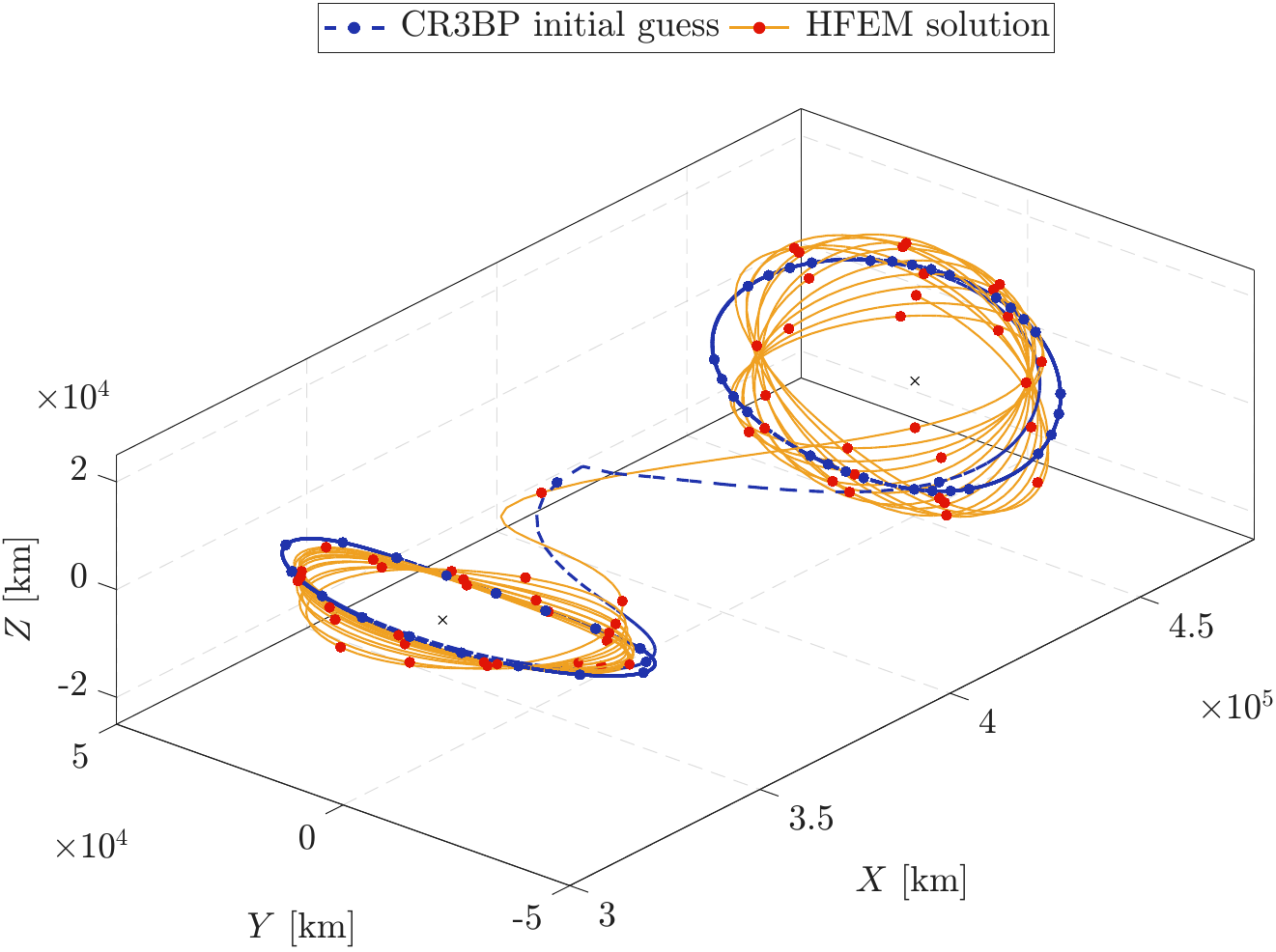}
			\caption{Converged LM solution considering no proximity constraints.}
			\label{fig:transfer_no_prox}
		\end{subfigure}
		\begin{subfigure}[c]{\linewidth}
			\centering
			\includegraphics[width=.83\textwidth]{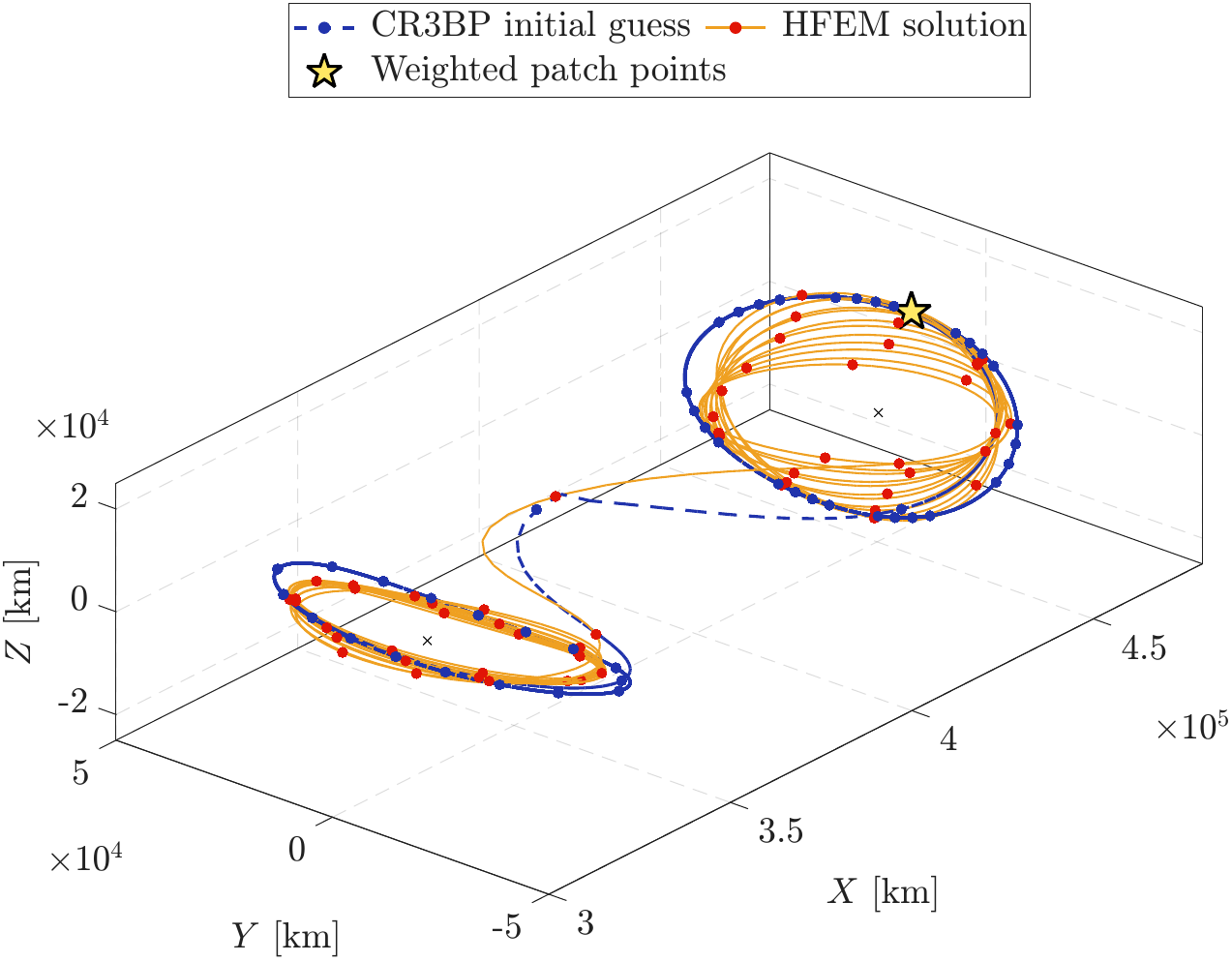}
			\caption{Converged LM solution weighting first error in position.}
			\label{fig:transfer_1st_prox}
		\end{subfigure}
	\end{minipage}
	\hfill
	\begin{minipage}[c]{.48\linewidth}
		\begin{subfigure}[c]{\linewidth}
			\centering
			\includegraphics[width=\textwidth]{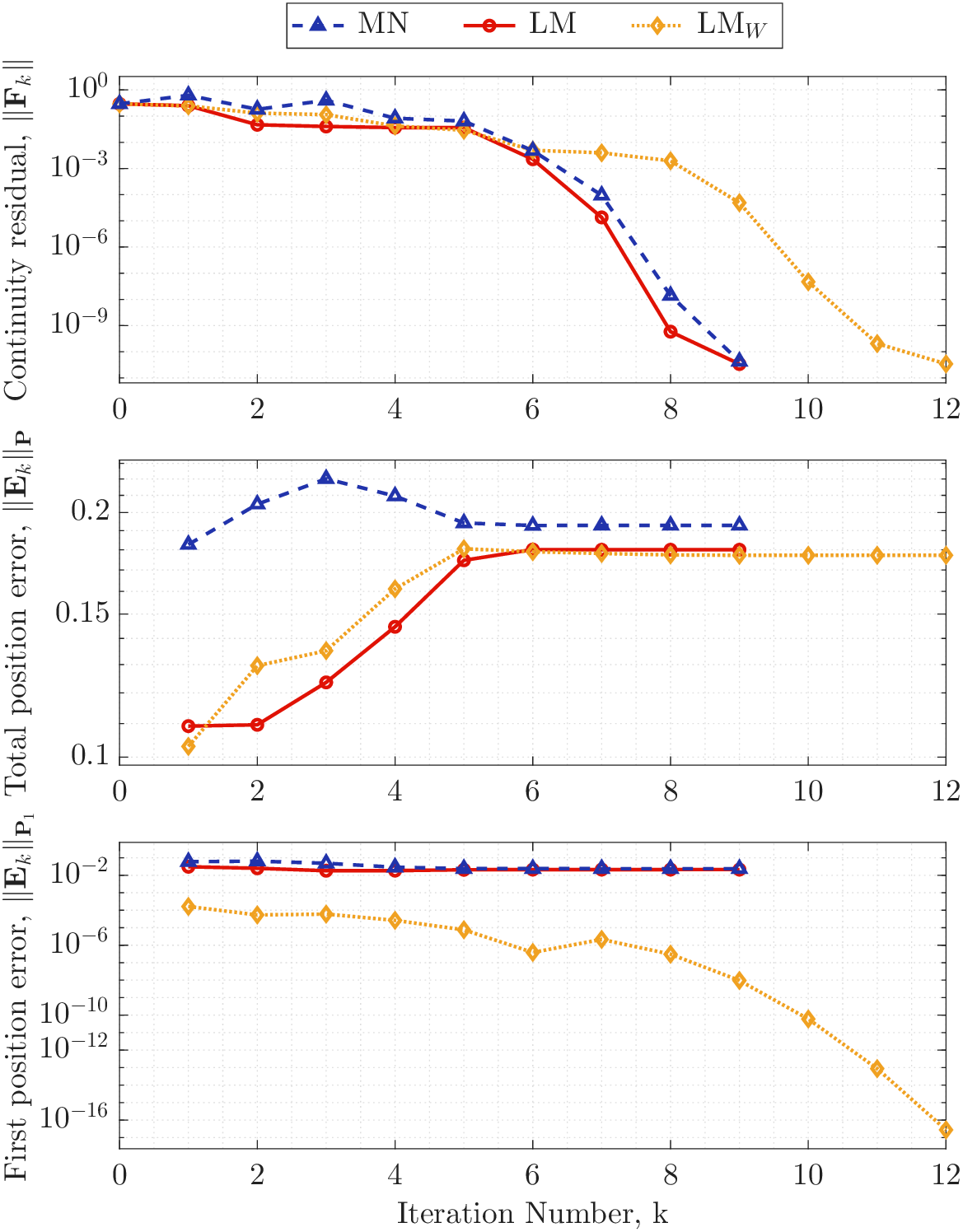}
			\caption{Evolution of meaningful parameters stemming from the application of the transition algorithms.}
			\label{fig:transfer_prox_study}
		\end{subfigure}
	\end{minipage}
	\caption{Application of the MN, LM, and weighted LM algorithms to transition the CR3BP transfer trajectory between the L2 Northern Halo orbit and L1 Lyapunov orbit to the HFEM considering sixty patch points. LM approaches consider an initial damping value of $\beta_0=10^{-5}$.}
	\label{fig:transfer_res}
\end{figure}

We consider the case presented in the second line of Table~\ref{tab:results_transfer} for further analysis. In this sense, Fig.~\ref{fig:transfer_no_prox} presents the converged solution stemming from the application of the LM algorithm considering 60 patch points, an initial damping of $\beta_0=10^{-5}$, and no proximity weighting. This solution is almost analogous to the result of the MN approach, which is hence omitted. By observation of the figure it is possible to verify that the converged solution has to deviate significantly from the initial guess in order to accommodate the challenging transfer between the two orbits of very different nature. In fact, we observe that having considered multiple revolutions around each orbit before and after the transfer serves the intended purpose of facilitating this maneuver, which is especially evident near the L1 Lyapunov orbit, where each successive revolution brings the spacecraft towards the $X$--$Y$ plane where the orbit resides. However, one possible concern that may arise from the observation of Fig.~\ref{fig:transfer_no_prox} is that the initial portion of the trajectory, near the L2 Northern Halo orbit, does not lie sufficiently close to the CR3BP solution, which might be relevant for a mission where a spacecraft is known to depart from this specific trajectory. This may be confirmed through the analysis of Fig.~\ref{fig:transfer_prox_study} where, besides the continuity residual and total position error, the evolution of the position error associated with the first patch point, denoted $\norm{\mathbf{E}_k}_{\mathbf{P}_1}$, is provided and compared with the MN alternative. To potentially address this concern, the introduction of the proximity constraint is finally studied, namely by taking $\mathbf{Q}=\text{diag}(1,1,1,0,\dots,0)=\mathbf{P}_1$ to penalize position deviations from the first patch point. The evolution of the aforementioned parameters under the proximity-weighted LM algorithm, denoted LM$_W$, is presented in Fig.~\ref{fig:transfer_prox_study} and the solution resulting from its application is provided in Fig.~\ref{fig:transfer_1st_prox}, with the first patch point highlighted by a star. Note that convergence was possible without the need for employing the adaptive weighting procedure detailed in Section~\ref{subsec:extensions}.

The analysis of Fig.~\ref{fig:transfer_res} highlights the difficulties in closely matching the two orbits, which differ so significantly from a spatial viewpoint. Furthermore, it exposes the advantage of the LM algorithm at controlling the shape of the converged solution through the inclusion of the proximity goal. By comparison of Fig.~\ref{fig:transfer_no_prox} and Fig.~\ref{fig:transfer_1st_prox} we verify that the proposed algorithm is able to converge to a solution that fully respects the additional objective by weighting the proximity to the initial patch point, matching the location of the departure point exactly at its corresponding position under the CR3BP guess. Fig.~\ref{fig:transfer_prox_study} attests to this result, highlighting that $\norm{\mathbf{E}_k}_{\mathbf{P}_1}$ is virtually brought to zero without any sacrifice to the total position error, which initially follows an evolution almost equal to the unconstrained case. Moreover, the number of iterations necessary to reach convergence is only slightly increased when accommodating the additional proximity goal, which was trivially imbued in the update step in a relaxed manner. Naturally, in this specific case where a single position constraint is being considered, one could have alternatively incorporated the proximity objective at the constraint level, as detailed in Section~\ref{subsec:MN}. This could level the grounds between the LM and MN methods, as both techniques allow for the inclusion of one constraint in three-dimensional position. However, we recall that such an approach would require further work on the construction of the optimization problem's variables and Jacobian. Additionally, it is less scalable, given the restrictions in degrees of freedom between $\mathbf{X}$ and $\mathbf{F}(\mathbf{X})$, as previously discussed. This is investigated in the next example.

\begin{figure}[h!]
	\centering
	\begin{minipage}[c]{\linewidth}
		\begin{subfigure}[c]{.48\linewidth}
			\centering
			\includegraphics[width=.9\textwidth]{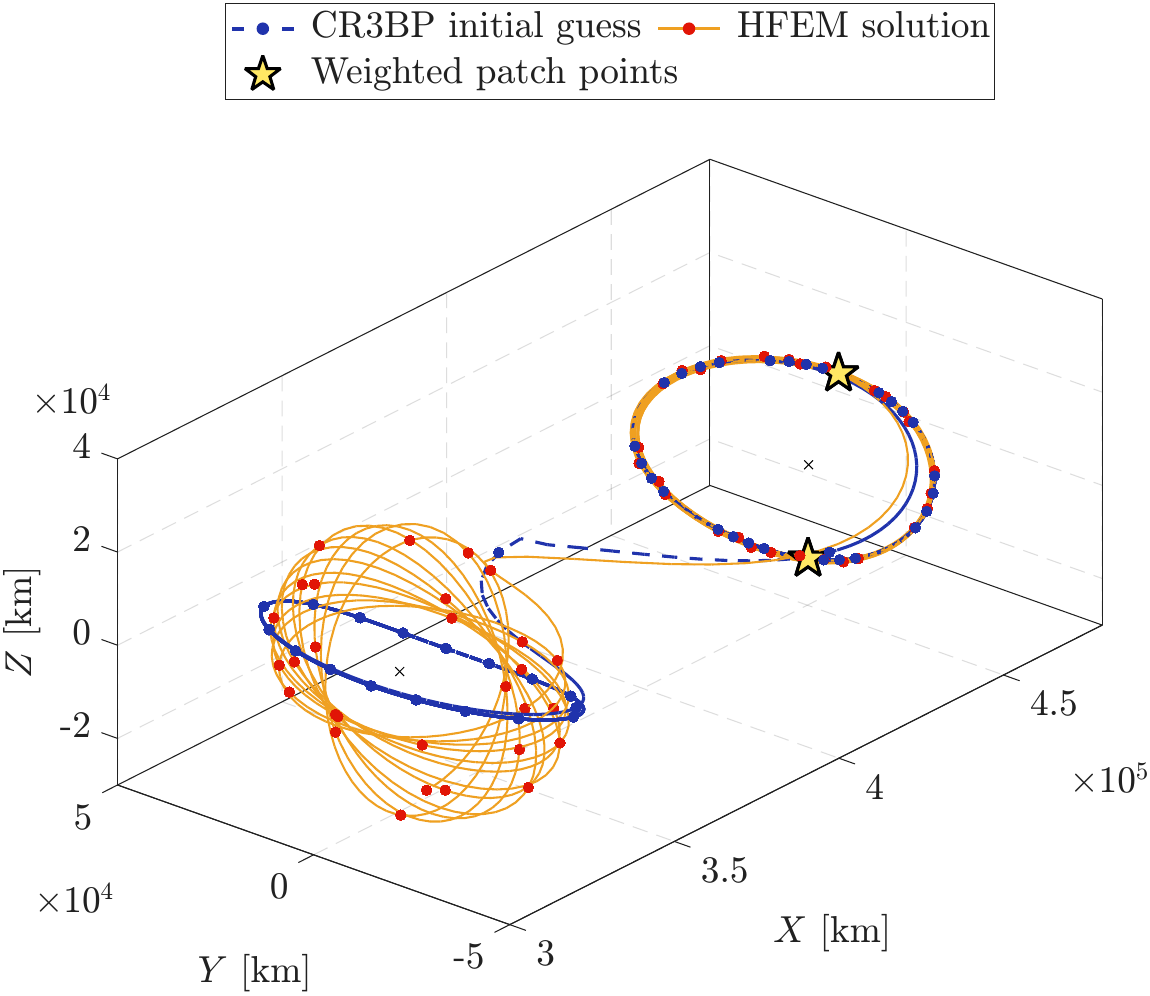}
			\caption{Converged LM solution weighting proximity to the L2 Northern Halo orbit.}
			\label{fig:LM_1_2_weight}
		\end{subfigure}
		\begin{subfigure}[c]{.48\linewidth}
			\centering
			\includegraphics[width=0.9\textwidth]{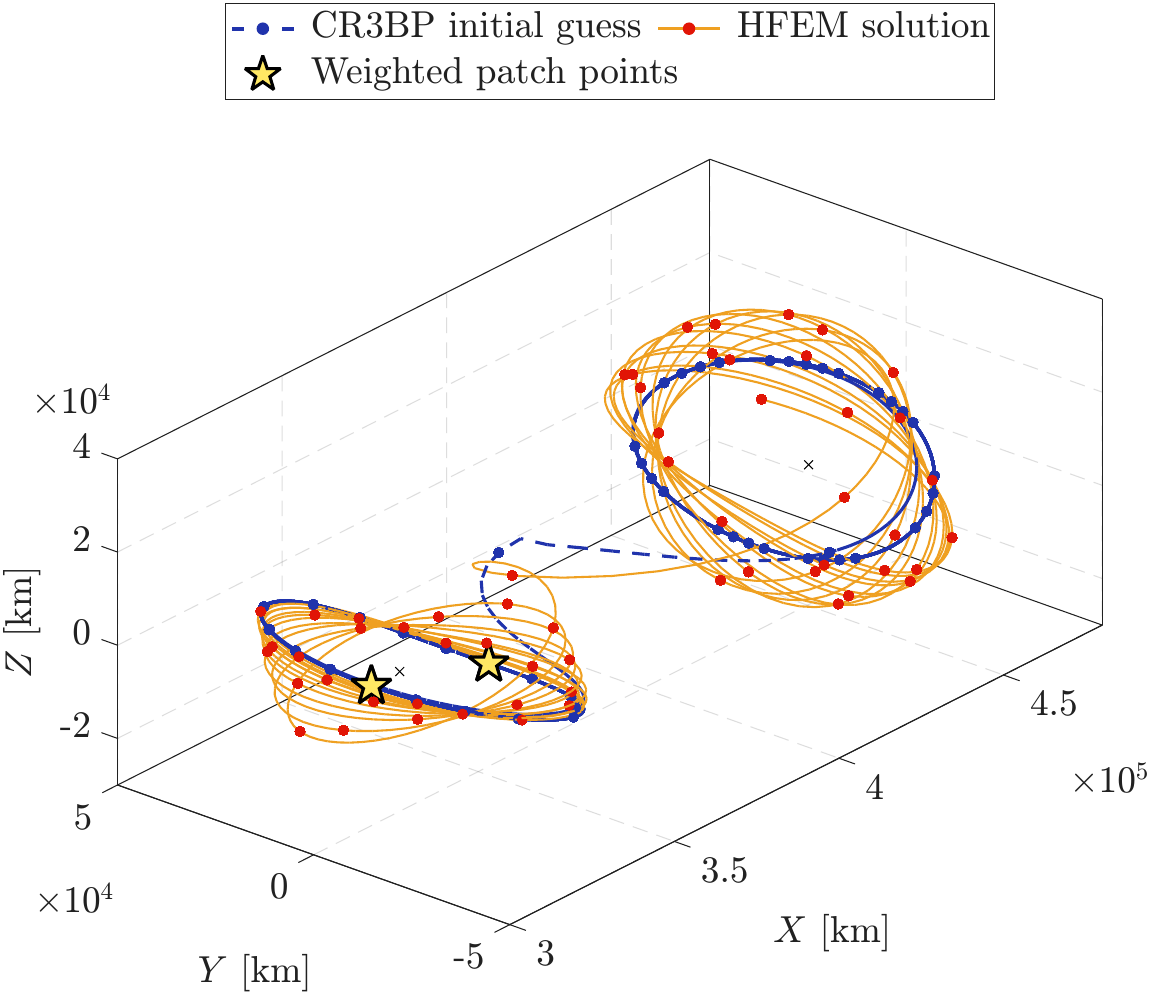}
			\caption{Converged LM solution weighting proximity to the L1 Lyapunov orbit.}
			\label{fig:LM_1_2_LAST_weight}
		\end{subfigure}
	\end{minipage}
	\hfill
	\begin{minipage}[c]{\linewidth}
		\begin{subfigure}[c]{\linewidth}
			\centering
			\includegraphics[width=0.9\textwidth]{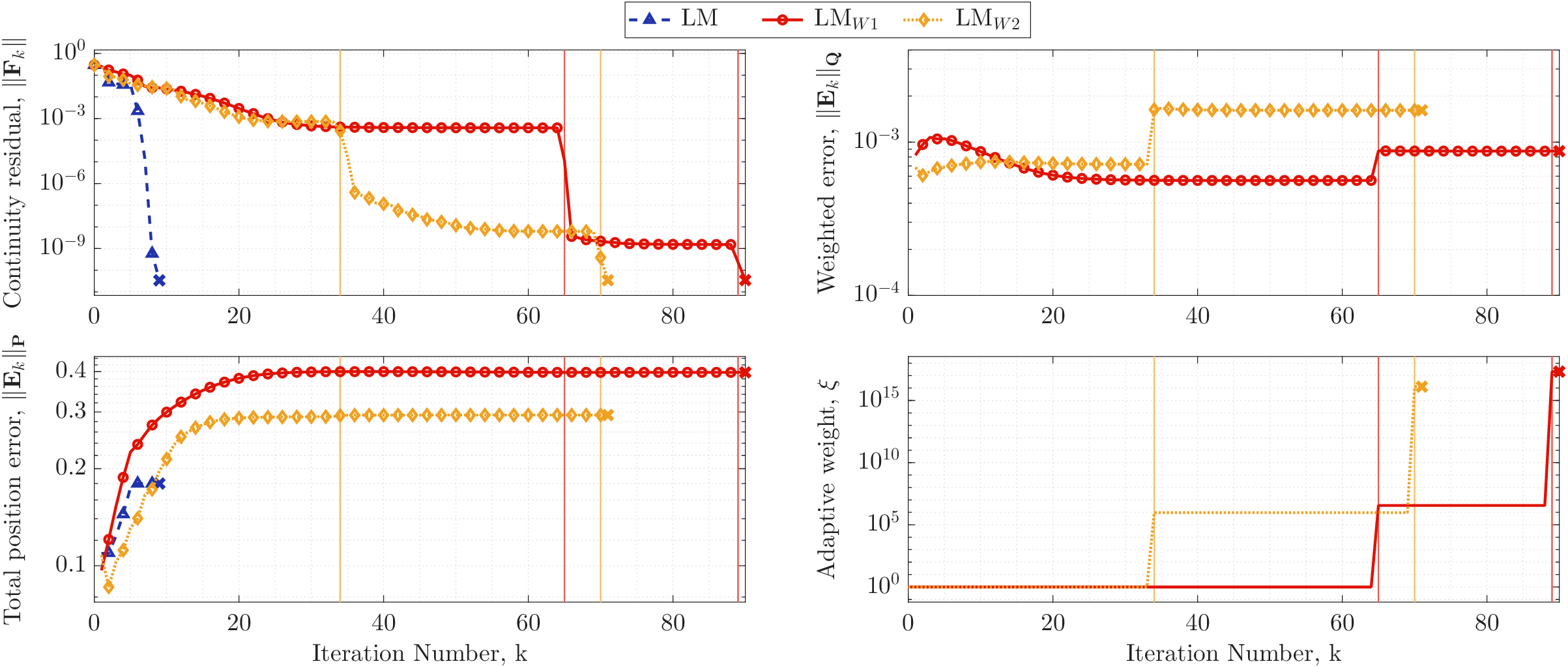}
			\caption{Evolution of meaningful parameters stemming from the application of the transition algorithms.}
			\label{fig:LM_comparison_weights}
		\end{subfigure}
	\end{minipage}
	\caption{Application of differently weighted LM algorithms to transition the CR3BP transfer trajectory between the L2 Northern Halo orbit and L1 Lyapunov orbit to the HFEM considering sixty patch points and different proximity objectives. All approaches consider an initial damping value of $\beta_0=10^{-5}$.}
\end{figure}

To further evaluate the capabilities of the proximity-constrained LM algorithm, it was pondered if it is possible to force convergence towards a solution that fully respects one of the orbits to a greater deal of accuracy, sacrificing proximity to the other. Recalling the previous example in Fig.~\ref{fig:transfer_1st_prox}, we note that respecting a single constraint in position along a desired orbit is insufficient to coerce the final converged trajectory into closely following it. Hence, to achieve the desired goal, it is necessary for the proximity objective to be imposed at multiple patch points along the CR3BP orbit to be tracked. In this case, a fixed-time formulation of the problem may not offer sufficient degrees of freedom to incorporate this objective at the constraint level. Thus, in order to proceed with using the MN update equation without running into feasibility issues, a reformulation of the multiple-shooting problem under a variable-time perspective would typically be necessary. On the contrary, the relaxed approach pursued by the proposed LM alternative ensures that this objective may be met, as best as possible, without abandoning the fixed-time formulation of the problem. In this sense, it was found that weighting the position entries of two consecutive points in a given orbit is sufficient to arrive to a final solution that respects it to a higher degree. Hence, two scenarios were considered: \textit{(i)} weighting of the first two patch points, from the L2 Northern Halo orbit, and \textit{(ii)} weighting of the last two patch points, from the L1 Lyapunov orbit. The corresponding weights in $\mathbf{Q}$ were set to $1$, while the remaining entries were kept null. To fulfill these requirements, employing the adaptive weighting approach detailed in Section~\ref{subsec:extensions} was found to be necessary. The adaptive weight $\xi$ is initially set to $1$ and then automatically adjusted at the beginning of each outer iteration, recalling Eq.~\eqref{eq:kappa}, until the algorithm is able to converge towards a solution that brings the continuity residual to $\norm{\mathbf{F}}<10^{-10}$. To reflect this change, the specified initial damping factor, $\beta_0$, is multiplied by $\xi$ at the start of each new outer iteration. When solving each optimization problem, the algorithm is stopped once there is a relative variation to the continuity residual between successive iterations lower than $10^{-4}$, 100 iterations have elapsed, or the continuity residual is brought below $10^{-10}$, in which case the final solution has been found. The results stemming from the application of the LM algorithm under the scenarios \textit{(i)} and \textit{(ii)}, which we denote LM$_{W1}$ and LM$_{W2}$, are provided in Fig.~\ref{fig:LM_1_2_weight} and Fig.~\ref{fig:LM_1_2_LAST_weight}, respectively. The corresponding evolution of the continuity residual, total position error, error of the weighted entries of the design vector, and adaptive weight are presented in Fig.~\ref{fig:LM_comparison_weights}, with markers at every second iteration for improved visibility. The evolution of the same parameters under the base LM algorithm with no proximity constraints, associated to the trajectory in Fig.~\ref{fig:transfer_no_prox}, is recalled for the sake of comparison. The corresponding curves are plotted until convergence of each algorithm is reached, which is evidenced by a cross. The vertical lines serve to indicate when a new outer iteration starts, i.e. once a new optimization problem, with adjusted weight $\xi$, is being solved.

The examples in Fig.~\ref{fig:LM_1_2_weight} and Fig.~\ref{fig:LM_1_2_LAST_weight} attest to the versatility of the algorithm proposed which, with a simple change to the weight matrix $\mathbf{Q}$, is able to converge to vastly different solutions that answer the desired objectives while retaining some measure of proximity to the initial guess. Furthermore, even when proximity to one of the orbits is imposed, the revolutions about the other orbit nonetheless remain in a contained region around it. This is confirmed in Fig.~\ref{fig:LM_comparison_weights} where, in spite of an increase to the total position error when compared to the unweighted LM baseline, the algorithms LM$_{W1}$ and LM$_{W2}$ converge to a solution that keeps $\norm{\mathbf{E}_k}$ within the same order of magnitude. The analysis of the continuity residual and weighted error exposes the counteracting contributions to the objective function of the optimization problem, i.e. as new optimization problems are formulated to more heavily penalize large continuity residuals, evidenced by the increase of $\xi$, the algorithms successively converge towards solutions that meet the proximity goals to a less desirable standard. Nonetheless, this gradual approach ensures that the final solution retains a measure of the two effects, at the cost of a considerably larger number of iterations necessary to reach convergence. Hence, while these examples are most relevant from an academic standpoint, they serve to illustrate the capabilities of the LM implementation proposed, namely its simplicity and versatility at handling proximity objectives without requiring a reformulation of the optimization problem.

\section{Concluding remarks}
\label{sec:conclusions}

This work proposes a novel algorithm for the transition of spacecraft trajectories within cislunar space from low-fidelity dynamical models, such as the CR3BP, to high-fidelity counterparts that consider the gravitational influence of the Earth, Moon, and Sun, whose trajectories are modeled via accurate ephemeris data. The proposed implementation pursues the use of a Levenberg--Marquardt algorithm, through the inclusion of an update step damping term in the formulation of the nonlinear least-squares optimization problem that stems from a fixed-time multiple-shooting viewpoint. This bestows the algorithm with a measure of inherent robustness and stability that proves useful when the initial guess is poor --- cases where the typical MN alternative from literature may diverge or converge to a trajectory farther away from the original CR3BP solution. The inclusion of additional constraints is pondered in a relaxed manner, namely to meet proximity objectives without hindering the integrity of the implementation, which is an advantage with respect to existing alternatives in the literature. An adaptive approach to the definition of the weights associated to the continuity and proximity related goals, at the objective function level, is covered to promote algorithmic stability and reduce manual adjustments.

Numerical results demonstrate the adequacy of the proposed LM algorithm when transitioning periodic orbits from the CR3BP orbit into quasi-periodic counterparts under the HFEM. The case of an L2 Northern Halo QPO demonstrates that LM is able to achieve similar performance to the MN alternative for higher-resolution initial guesses, but retains better proximity to the original orbit when the initial guess is poor, through an appropriate selection of the initial damping factor. These results were generalized to various QPOs from different families near the L2 and L1 Lagrange points. Doing so further evidences that the MN baseline encounters numerical difficulties in certain cases, leading to divergence if the update step magnitude is not limited. Even then, it was found that the typical brute-force attenuation pursued in the literature may be unpredictable or significantly slow down convergence speed, being objectively inferior to the automatic adjustment of the damping factor of the proposed LM algorithm. When transitioning a transfer trajectory between two periodic orbits from the CR3BP, similar conclusions in terms of adequacy and robustness were drawn. Furthermore, the benefits associated with the ease of including proximity constraints in the LM formulation were demonstrated, allowing for the converged solution to be shaped in accordance with various objectives without requiring a reformulation of the optimization problem. The adaptive weight definition was found to be crucial in the cases where the LM algorithm initially converged to a poorly-continuous solution, by increasingly penalizing the continuity residual. Employing the adaptive procedure was however verified to significantly slow convergence speed. Still, the robustness and versatility of the proposed LM package, all while remaining simple to implement and computationally light to execute, make it a valuable alternative to the MN update scheme for practitioners focusing on HFEM transition.

Future work should assess the behavior the proposed LM algorithm under a variable-time formulation of the multiple-shooting problem, comparing it with the fixed-time results presented in this paper and the MN baseline. The use of alternative formulations of the optimization problem besides the unconstrained least-squares, such as QCQP, should also be investigated. Furthermore, the application of the proposed algorithm should be extended to orbits whose period resonates with the dynamics or whose periapsis is small (e.g. near-rectilinear Halo orbits), which are notoriously challenging to transition \cite{park2025CharacterizationL2Analogs,Sanaga2024ChallengingRegion} --- where the LM scheme may encounter difficulties due to its inherently monotonous nature. This is expected to require the use of additional techniques, such as homotopy, to gradually bridge the CR3BP and HFEM dynamics \cite{sanaga2026UTSTori,soto2023thesis}.

\section*{Acknowledgments}
This work was supported by LARSyS FCT funding \linebreak(UID/50009/2025: DOI 10.54499/UID/50009/2025; LA/P/0083/2020: DOI 10.54499/LA/P/0083/2020), and by NEURASPACE Project, Contract No. 9, under Regulation (EU) 2021/241 of the European Parliament and of the Council of February 12, 2021 and the Portuguese Recovery and Resilience Program (PRR), in component 05 - Capitalization and Business Innovation, under Notice No. 425 01/C05-i01/2021 of the Regulation of Mobilizing Agendas/Alliances for re-industrialization.

\appendix
\section{Parametric study}
\label{ap:param_study}
The LM algorithm proposed in this work, in Section~\ref{subsec:LM}, considers two factors, $0<\alpha<1$ and $\eta>1$, that are used to adjust the damping term, $\beta_k$, according to Algorithm~\ref{alg:LM}. These parameters are, for the most part, defined in an empirical or heuristic manner, and may require prior tuning for each application. In the literature, various strategies can be found, with some authors preferring to scale $\beta_k$ according to reciprocal factors, say $(\alpha,\eta)=(0.1,10)$ or $(\alpha,\eta)=(0.5,2)$ for a more extreme and moderate choice, respectively. Another standard alternative, termed \textit{delayed gratification} \cite{Transtrum2011LM}, is to take $\alpha<1/\eta$, with the intent of rewarding successful steps and adjusting cautiously when the step fails. Example choices for delayed gratification are $(\alpha,\eta)=(0.33,2)$ and $(\alpha,\eta)=(0.2,1.5)$.

\begin{figure}[h!]
	\centering
	\includegraphics[width=0.97\textwidth]{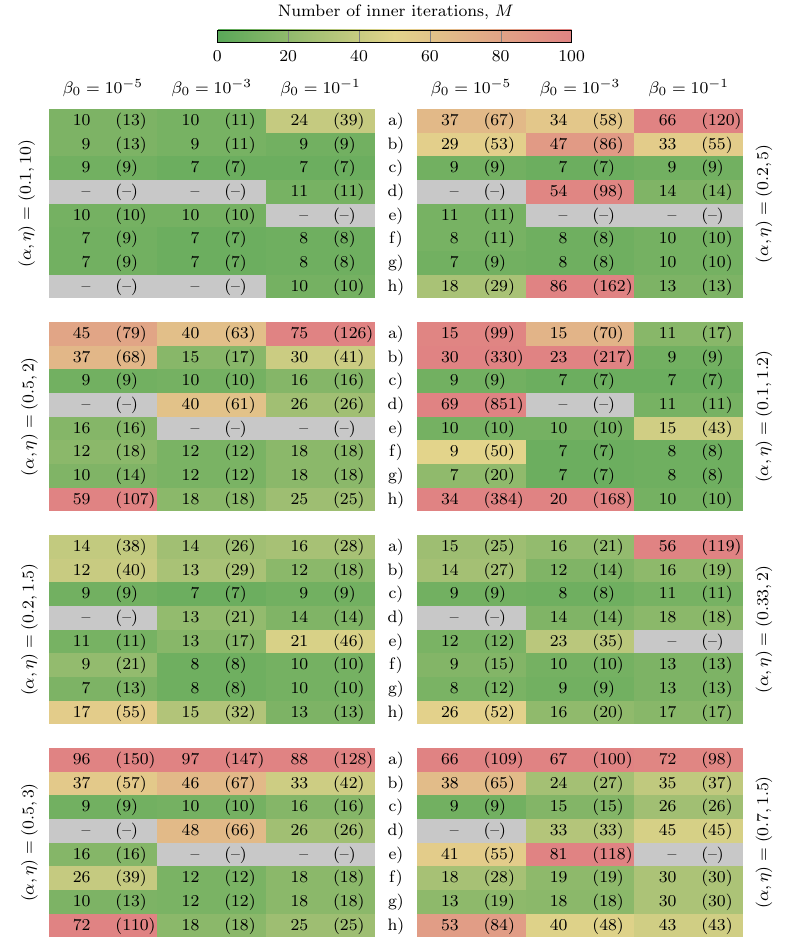}
	\caption{Performance map for the total number of outer (and inner iterations), $m~(M)$, to reach convergence when transitioning the CR3BP orbits of Fig.~\ref{fig:more_orbits} with the proposed LM algorithm, for various values of the initial damping, $\beta_0$, and its adjustment factors $(\alpha,\eta)$. Cell color varies according to the number of inner iterations, $M$.}
	\label{fig:parametric_study_map_N}
\end{figure}

\begin{figure}[h!]
	\centering
	\includegraphics[width=0.97\textwidth]{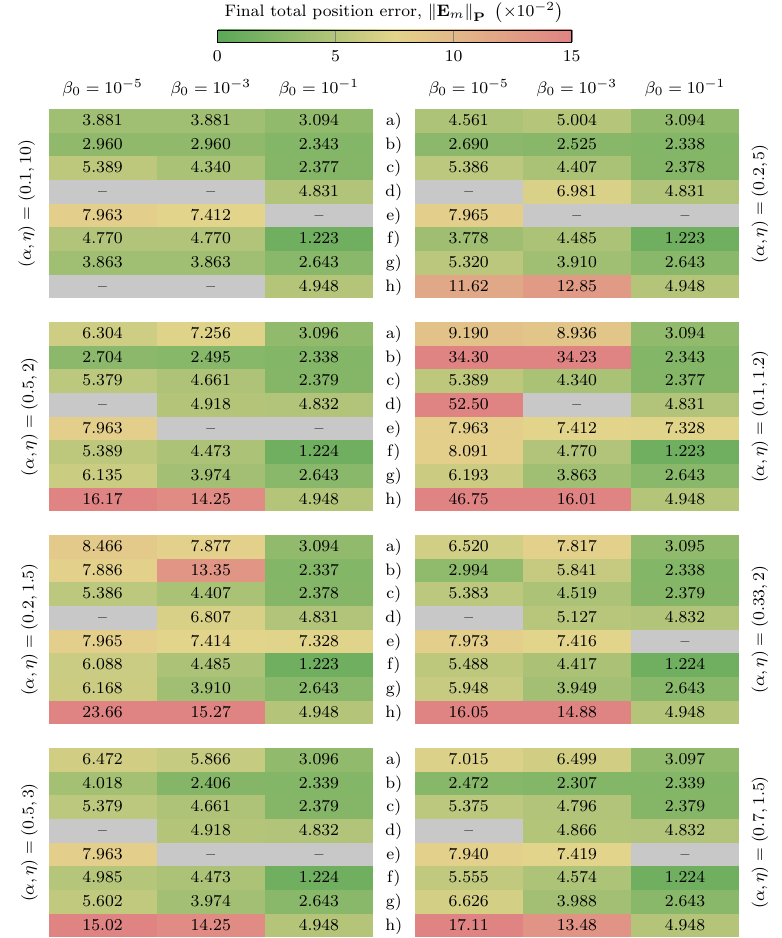}
	\caption{Performance map for the total position error of the converged solution, $\norm{\mathbf{E}_m}_\mathbf{P}$, when transitioning the CR3BP orbits of Fig.~\ref{fig:more_orbits} with the proposed LM algorithm, for various values of the initial damping, $\beta_0$, and its adjustment factors $(\alpha,\eta)$.}
	\label{fig:parametric_study_map_E}
\end{figure}

In order to determine an adequate $(\alpha,\eta)$ choice for the purposes of HFEM transition in this work, a preliminary analysis was carried out considering the transition of the orbits from Table~\ref{tab:more_orbits_periods}. Following common strategies from literature for the definition of $(\alpha,\eta)$, the reciprocal pairs $(0.1,10)$, $(0.2,5)$, and $(0.5,2)$ are studied, as well as the pairs $(0.1,1.2)$, $(0.2,1.5)$, and $(0.33,2)$ for a delayed gratification approach. In addition, with the intent of evaluating an opposing alternative to delayed gratification, the pairs $(0.5,3)$ and $(0.7,1.5)$ were also evaluated for representing choices where $\alpha>1/\eta$. Fig.~\ref{fig:parametric_study_map_N} presents a performance map, for each of these $(\alpha,\eta)$ pairs, based on the number of total outer iterations, $m$, and inner iterations, $M$, in parenthesis. For each $(\alpha,\eta)$ pair, the initial damping factor is selected from $\beta_0=(10^{-5},10^{-3},10^{-1})$. The color of each cell is defined according to the color bar atop the figure, in terms of $M$, since it provides a better overall picture of the total computational cost. In Fig.~\ref{fig:parametric_study_map_E}, a similar performance map is provided for the total position error of the corresponding converged solutions, $\norm{\mathbf{E}_m}_\mathbf{P}$. In case convergence is not reached within $m=100$ outer iterations, the respective value is removed from both maps.

The analysis of Fig.~\ref{fig:parametric_study_map_N} and Fig.~\ref{fig:parametric_study_map_E} evidences that the choice of $\alpha$ and $\eta$ bears an important impact on how quickly the proposed algorithm converges to a final solution. In this sense, while these parameters do not impair the robustness of the proposed scheme \textit{per se}, a poor choice may significantly slow down convergence --- beyond the 100 outer iteration limit imposed here. In particular, one may observe that taking $(\alpha,\eta)=(0.1,10)$ may prove unreliable. Despite providing, on average, the fastest responses and converging to HFEM trajectories that are closest to the initial CR3BP guess, it is also the case that this choice results in the highest amount of runs where the algorithm is unable to converge within $100$ iterations. Further analysis shows that this occurs due to the algorithm stalling, which is potentiated by the fact that $\alpha$ and $\eta$ are reciprocal. Such choices, especially if $\alpha$ is small and $\eta$ is large, motivate $\beta_k$ to remain constant unless changing is absolutely necessary. Hence, the algorithm is not incentivized to adjust (namely, to decrease) the damping factor, leading to extremely slow convergence. The same reasoning may be invoked to justify the lack of convergence within 100 iterations, for three test cases, when $(\alpha,\eta)=(0.2,5)$ and $(0.5,2)$. On the contrary, the choices corresponding to a delayed gratification approach --- i.e., $(\alpha,\eta)=(0.1,1.2)$, $(0.2,1.5)$, and $(0.33,2)$ --- are usually better at ensuring convergence. In fact, only for $(\alpha,\eta)=(0.1,1.2)$ was it possible to achieve convergence for orbit d), with $\beta_0=10^{-5}$, though at a significant cost in terms of total inner iterations and final position error, which is the largest of all test-cases approached\footnote{The generalized difficulty in convergence for orbit d) with $\beta_0=10^{-5}$, observed for every $(\alpha,\eta)$ option, suggests a particular lack of quality of the linear approximation at play. As analyzed in Section~\ref{subsec:more_orbits}, this prompts the algorithm to inaccurately divert significantly from the initial CR3BP guess. Fortunately, this issue is easily resolved with a larger choice of $\beta_0$, i.e. by reducing the level of trust initially put on the linear approximation.}. Despite this, $(\alpha,\eta)=(0.2,1.5)$ and $(0.33,2)$ are still, arguably, the better choices for delayed gratification, since they retain lower computational costs --- nearly on par with the fastest option, $(\alpha,\eta)=(0.1,10)$. Observation of the total final position error, however, suggests that the delayed gratification approach typically leads to HFEM trajectories that are farther away from the initial CR3BP guess. This is a direct consequence of incentivizing successful steps which, whenever possible, will tend to decrease $\beta_k$ and inevitably divert the algorithm towards a solution that is farther away. Finally, the options of $(\alpha,\eta)=(0.5,3)$ and $(0.7,1.5)$, selected to contrast with delayed gratification, were found to slightly improve the final position errors, though at a significant increase in terms of numerical cost, especially evident for the case of orbit a).

\begin{figure}[h]
	\centering
	\includegraphics[width=0.98\textwidth]{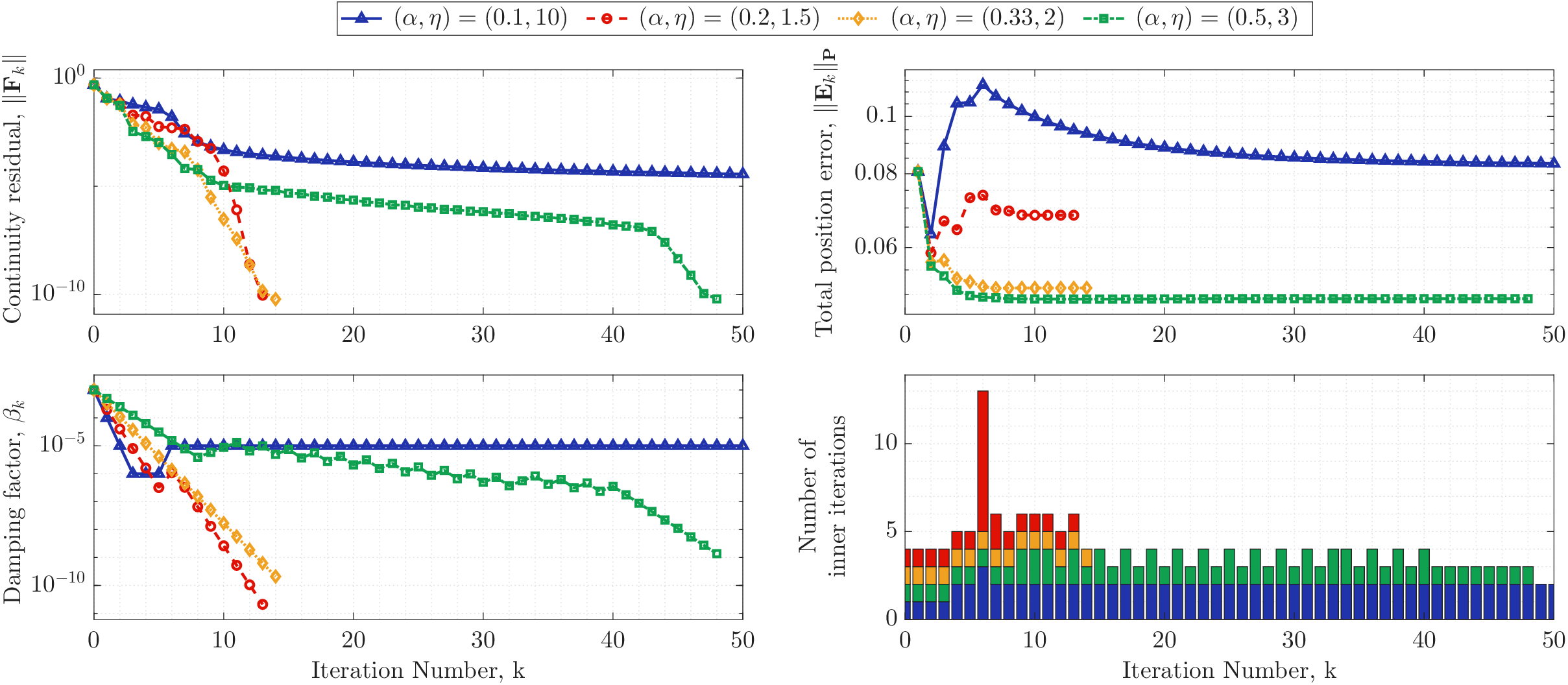}
	\caption{Number of inner LM iterations and evolution of the damping factor for orbit d) with $\beta_0=10^{-3}$ and various choices of $(\alpha,\eta)$, over the first 50 outer iterations.}
	\label{fig:parametric_study_ex}
\end{figure}

Given the discussion above, it can thus be said that the choice of the specific approach pursued to select $\alpha$ and $\eta$ ultimately amounts to a trade-off between algorithmic reliability and the quality of the final solution. To further illustrate these considerations, Fig.~\ref{fig:parametric_study_ex} presents the evolution of the continuity residual, total position error, damping factor, and number of inner iterations at each outer iteration, for orbit d) with $\beta_0=10^{-3}$. The figure focuses on the choices of $(\alpha,\eta)=(0.1,10)$, $(0.2,1.5)$, $(0.33,2)$, and $(0.5,3)$, and is limited to the first 50 outer iterations, for clarity.

The observation of Fig.~\ref{fig:parametric_study_ex} exposes how the internal updates to the damping factor may influence the LM algorithm's behavior. In this example, taking both $(\alpha,\eta)=(0.1,10)$ and $(0.5,3)$ leads to significantly slower convergence, although for two distinct reasons. In the case of the former, as previously discussed, the algorithm stalls and the damping factor remains constant, not allowing for larger update steps to be applied. We verify that the algorithm may be put in this undesirable mode due to an overreaching trust on the linear approximation during the initial iterations, which diverts it from the original CR3BP orbit. This is ultimately a consequence of the small value of $\alpha$, which promotes an accentuated decrease in $\beta_k$ at the start. On the contrary, $(\alpha,\eta)=(0.5,3)$ directs the algorithm to a solution that is closer to the initial CR3BP guess, by restricting changes to the damping factor, though the inertia to decrease $\beta_k$ is seen to come at the cost of a slow convergence. For these reasons, we verify that the choice with $\alpha=0.3$ proves the most adequate in this specific application case, as it leads to a monotonous behavior in terms of $\beta_k$. This potentiates a quick convergence towards a close-by solution that is almost on par with the choice of $\alpha=0.5$, in terms of $\norm{\mathbf{E}_m}_\mathbf{P}$. In comparison, $(\alpha,\eta)=(0.2,1.5)$ tracks a close second, requiring only additional inner iterations in two of its outer steps, but converging to a slightly farther solution. In fact, as suggested by the previous data in Fig.~\ref{fig:parametric_study_map_N} and Fig.~\ref{fig:parametric_study_map_E}, the choices of $(\alpha,\eta)=(0.33,2)$ and $(0.2,1.5)$ lead to a behavior that is similar to the one observed here: the former typically converges to a solution that is closer to the initial CR3BP guess, but the latter benefits from its quick adaptability, generally presenting a higher rate of success and converging in fewer iterations (contrary to the specific case in Fig.~\ref{fig:parametric_study_ex}). In this sense, one could argue in favor of any of these two options, depending on what strategy best fits the specific needs of the application. In this work, both options are considered, as detailed throughout the paper. Naturally, more complex procedures to the definition of $\beta_k$ could potentially be explored in the future, as well as other $(\alpha,\eta)$ factors, though this analysis was deemed outside the scope of the present investigation.

\bibliographystyle{elsarticle-num} 
\bibliography{refs.bib}
\end{document}